\documentclass[preprint,12pt]{elsarticle}

\usepackage[latin1]{inputenc}
\usepackage{setspace}
\usepackage{amsmath}
\usepackage{amsfonts}
\usepackage{amssymb}
\usepackage{amsthm}
\usepackage{amssymb}
\usepackage{graphicx}
\usepackage{hyperref}
\usepackage{multirow}
\usepackage{longtable}
\usepackage{resizegather}
\usepackage[dvipsnames]{xcolor}
\usepackage[normalem]{ulem}
\usepackage{cancel}
\usepackage{lineno}
\usepackage{cases}
\usepackage{tabularx,longtable,booktabs,multicol}

\newtheorem{theorem}{Theorem}
\newtheorem{corollary}[theorem]{Corollary}
\newtheorem{proposition}[theorem]{Proposition}
\newtheorem{lemma}[theorem]{Lemma}
\newtheorem{defin}[theorem]{Definition}

\newtheorem{remk}[theorem]{Remark}

\DeclareMathOperator{\tr}{tr}
\DeclareMathOperator{\ev}{ev}

\newcommand{\thesisonly}[1]{}
\newcommand{\omitthesisonly}[1]{}

\newcommand{\newdraft}[1]{}
\newcommand{\omitnewdraft}[1]{}

\newcommand{\omitarticleonly}[1]{}

\journal{ }

\begin{document}

\begin{frontmatter}

%% Title, authors and addresses

%% use the tnoteref command within \title for footnotes;
%% use the tnotetext command for theassociated footnote;
%% use the fnref command within \author or \address for footnotes;
%% use the fntext command for theassociated footnote;
%% use the corref command within \author for corresponding author footnotes;
%% use the cortext command for theassociated footnote;
%% use the ead command for the email address,
%% and the form \ead[url] for the home page:
%% \title{Title\tnoteref{label1}}
%% \tnotetext[label1]{}
%%\author{Name\corref{cor1}\fnref{label2}}
%%\ead{email address}
%\ead[url]{home page}
%%\fntext[label2]{}
%%\cortext[cor1]{}
%%\affiliation{organization={},
%%             addressline={},
%%             city={},
%%             postcode={},
%%             state={},
%%             country={}}
%%\fntext[label3]{}

%\title{On the $k$-independence number of a graph}
\title{The optimal bound on the $3$-independence number obtainable from a polynomial-type method}
%% use optional labels to link authors explicitly to addresses:
%% \author[label1,label2]{}
%% \affiliation[label1]{organization={},
%%             addressline={},
%%             city={},
%%             postcode={},
%%             state={},
%%             country={}}
%%
%% \affiliation[label2]{organization={},
%%             addressline={},
%%             city={},
%%             postcode={},
%%             state={},
%%             country={}}
%\author{Lord C. Kavi, Michael W. Newman}
%\author{Lord C. Kavi}
\author[inst1]{Lord C. Kavi\corref{cor1}}
\ead{lkavi060@uottawa.ca}
\cortext[cor1]{Corresponding author}
\affiliation[inst1]{organization={Department of Mathematics, University of Ottawa},%Department and Organization
            %addressline={Address One}, 
            city={Ottawa},
            %postcode={00000}, 
            state={ON},
            country={Canada}}

\author[inst1]{Mike Newman}
\ead{mnewman@uottawa.ca}
%\author[inst1,inst2]{Author Three}

%\affiliation[inst2]{organization={Department Two},%Department and Organization
            %addressline={Address Two}, 
            %city={City Two},
            %postcode={22222}, 
            %state={State Two},
            %country={Country Two}}

\begin{abstract}
%% Text of abstract
A $k$-independent set in a connected graph is a set of vertices such that any two vertices in the set are at distance greater than $k$ in the graph. The {$k$-independence number} of a graph, denoted $\alpha_k$, is the size of a largest $k$-independent set in the graph.	Recent results have made use of polynomials that depend on the spectrum of the graph to bound the $k$-independence number. They are optimized for the cases $k=1,2$. There are polynomials that give good (and sometimes) optimal results for general $k$, including case $k=3$. In this paper, we provide the best possible bound that can be obtained by choosing a polynomial for case $k=3$ and apply this bound to well-known families of graphs including the Hamming graph.
\end{abstract}
%Graph; $k$-Independence number; 	Spectrum;	Interlacing; Polynomials.

\begin{keyword}
%% keywords here, in the form: keyword \sep keyword
Graph \sep $k$-Independence number \sep Spectrum \sep Interlacing \sep Hamming graph
%% PACS codes here, in the form: \PACS code \sep code
%\PACS 0000 \sep 1111
%% MSC codes here, in the form: \MSC code \sep code
%% or \MSC[2008] code \sep code (2000 is the default)
\MSC[2010] 05C50 \sep 05C69
\end{keyword}

\end{frontmatter}

%% \linenumbers

%% main text

\section {Introduction}
An {independent set}, also known as  a stable set or coclique, in a graph is a set of vertices, no two of which are adjacent. The size of a largest independent set is called the {independence number}. Two classical eigenvalue bounds on the independence number are the Hoffman's ratio bound~\cite{Haemers} and Cvetkovi\'{c}'s inertia bound~\cite{Cvetkovic}. There are some generalizations of the notion of independent set of a graph; of interest to us is the following: 
A {$k$-independent set} in a graph is a set of vertices such that any two vertices in the set are at a distance of at least $k+1$ in the graph. The {$k$-independence number} of a graph, denoted $\alpha_k$, is the size of a largest $k$-independent set in the graph. %set of vertices such that any two vertices in the set are at distance at least $k+1$. the k-independence number of a graph, which is the maximum size of a set of vertices at pairwise distance greater than k. the size of the largest set of vertices such that any two vertices in the set are at distance larger than k.

We remark that the $k$-independence number is related to coding theory (for instance, codes and anticodes are $k$-independent sets in the Hamming graph, (MacWilliams and Sloane~\cite{sloane})), the average
distance (Firby and Haviland~\cite{Firby}), packing chromatic number (Goddard, Hedetniemi, Hedetniemi, Harris, Rall~\cite{Goddard}), injective chromatic number (Hahn, Kratochv\'{i}l, \v{S}ir\'{a}\v{n}, Sotteau~\cite{Hahn}), and strong chromatic index  (Mahdian~\cite{Mahdian}) of a graph, as stated in Abiad, Cioab\u{a}, Tait \cite{Abiad1} and Fiol \cite{fiol2}. In particular, upper bounds on the $k$-independence number directly give lower bounds on the corresponding distance or packing chromatic
number (Abiad, Cioab\u{a}, Tait~\cite{Abiad1}).

Since computing $k$-independence number is NP-hard, see Kong and Zhao\cite{Kong}, it makes sense to look at spectral bounds instead.  Abiad, Cioab\u{a} and Tait~\cite{ Abiad1} gave inertial and ratio type spectral bounds on the  $k$-independence number. These were later generalized by Abiad, Coutinho, Fiol \cite{abiad2019k} using interlacing, which involves taking polynomial of degree at most $k$. In general any degree $k$ polynomial gives a bound on the $k$-independence number. A good bound therefore depends on making the right choice of a polynomial.  For the generalized Hoffman (or ratio type) bound for $\alpha_k$, the polynomial $p(x)=x$ gives the standard Hoffman bound for the independence number $\alpha_1$. Abiad, Coutinho and Fiol ~\cite{abiad2019k} also gave the right choice of polynomial for $k=2$, and proposed a polynomial for a general $k$. This polynomial, however, is often not the best choice. Fiol~\cite{fiol2} introduced the so-called minor polynomials for the same purpose. See also, Abiad, Coutinho, Fiol, Nogueira, and Zeijlemaker~\cite{Abiad2} where mixed integer linear programs (MILP) were used to optimize the choice of polynomials for the inertia type bound.  Abiad, Dalf\'{o}, Fiol and Zeijlemaker~\cite{Abiad33} also studied the relationship between the inertia and the ratio type bounds for $\alpha_k$. 

In this article, we provide the best possible polynomial to bound the $3$-independence number. We considered the polynomial $p(x)=x^3+bx^2+cx$, $b,c\in \mathbb{R}$. %To determine which eigenvalue(s) minimize the polynomial, 
For $b$ fixed, we give the eigenvalue that minimizes $p(x)$ when $c$ is within a given interval. This is explored in Lemma~\ref{lemma5}. With this polynomial, we are able to bound the $3$-independence number when $b$ is fixed, and then by the help of Lemma~\ref{lemma5}, we determine the optimal $c$ in terms of $b$. The bound is now in terms of $b$, and optimizing this bound enables us to obtain the optimal $b$ and $c$ for $p(x)$ in terms of the eigenvalues of the graph. This approach and the resulting optimal polynomial alongside its corresponding bound are in Theorem~\ref{Gmaink3}.  We compute the obtained bound explicitly for several graph families and compare them to previous results. For instance, we investigate tightness of this bound on the Hamming graph $H(d,q)$. In particular, we give a construction of $3$-independent sets in $H(d,2)$ and show tightness of the bound for $d=2^r$ and $d=2^r-1$ with $r\in \mathbb Z^+$.

In principle, our approach might extend to values of $k$ greater than 3  but would certainly become increasingly technical. It is not clear whether it is possible to explicitly  write down an optimal polynomial in terms of $k$.

\section{Background}
Let $G=(V, E)$ be a connected graph with $n =|V|$ vertices,  and adjacency matrix $A$ with the eigenvalues $\ev(G) =\{\theta_{0}, \theta_{1}^{m_1},\dots, \theta_{d}^{m_d}\}$, where the different eigenvalues are in decreasing order, $\theta_0> \theta_1>\dots >\theta_d$, and the superscripts stand for their multiplicities. When the eigenvalues are presented with possible repetitions, we shall indicate them by   $\ev(G) :\lambda_1\ge \lambda_2\ge \dots\ge \lambda_n$. 

Eigenvalue interlacing is an important tool used in proving most of the results in this paper. 	Let $A$ and $B$ be square matrices  with respective eigenvalues $\mu_1\ge \mu_2\ge \dots\ge \mu_k$ and $\lambda_1\ge \lambda_2\ge \dots\ge \lambda_n$, where $k< n$. We say the first sequence \textit{interlaces} the second sequence if $\lambda_i\ge \mu_i\ge \lambda_{n-k+i}$ for $i=1,\dots, k$. Moreover, the interlacing is said to be \textit{tight}, if for some $\ell$, with $1\le \ell \le l$, we have $\lambda_i=\mu_i$ for $i=1,\dots,l$ and $\mu_i=\lambda_{n-k+i}$ for $i=l+1,\dots, k$.
%We note that in Theorem~\ref{Pinterlacing}, matrix $R$ satisfies $R^TR=I$, and it turns out that this property is sufficient. The next result tells us more about the interlacing of the eigenvalues of $M$ by the eigenvalues of a matrix $N=R^TMR$. 

The following interlacing theorem is attributed to Haemers \cite{Haemers}. See also Godsil and Royle \cite{Godsil1}.
\begin{theorem}\label{Interlacing}(Interlacing \cite{Haemers, Godsil1})
	Let $A$ be a real symmetric $n \times n$ matrix,  $R$ a real $n \times k$ matrix such that $R^TR = I$, and $B = R^TAR$ for $n\ge k$. Let the eigenvalues of $A$ be $\lambda_1\ge \lambda_2\ge \dots \ge \lambda_n$ with corresponding orthonormal eigenvectors $\mathbf{v}_1 ,\mathbf{v}_2,\dots,\mathbf{v}_{n}$, and let the eigenvalues
	of $B$ be $\mu_1\ge \mu_2\ge \dots\ge \mu_k$  with corresponding orthonormal eigenvectors $\mathbf{u}_1 ,\mathbf{u}_2,\dots,\mathbf{u}_{k}$. Then the
	following statements hold.
	\begin{itemize}
		\item[i.] The eigenvalues of $B$ interlace the eigenvalues of $A$.
		\item [ii.] If $\mu_i=\lambda_i $ (or $\mu_i=\lambda_{n-k+i}$), then there exists $\mathbf{u}\in \mathbb{R}^k$ with $B\mathbf{u} = \mu_i\mathbf{u}$ and $A(R\mathbf{u}) = \mu_i(R\mathbf{u})$.
		\item[iii.] If for some $\ell$ we have $\mu_i=\lambda_i$ for $i=1,\dots, \ell$ (or $\mu_i=\lambda_{n-k+i}$ for $i=\ell,\dots, k$), then $A(R\mathbf{u}_i) =\mu_i(R\mathbf{u}_i)$ for  $i=1,\dots, \ell$ ($i=\ell,\dots, k$ respectively).
		\item[iv.] If the interlacing is tight, then $AR=RB$.
	\end{itemize}
\end{theorem}

Two interesting particular cases where interlacing occurs (obtained by choosing appropriately the matrix $R$) are the following.  Let $A$ be the adjacency matrix of a graph $G=(V,E)$. First, if $B$ is a principal submatrix of $A$, then $B$ corresponds to the adjacency matrix of an induced subgraph $G'$ of $G$. Second, given a partition $\pi$ of the vertices as $V=V_1\cup\dots\cup V_n$, the quotient matrix $B$ of $A$ with respect to partition $\pi$ is the matrix with elements $b_{ij}$, for $i,j=1,\dots,n$, being the average row sums of the corresponding block $A_{ij}$ of $A$. The quotient matrix $B$ is similar to $R^TAR$. Moreover, if the interlacing is tight, Theorem \ref{Interlacing}$(iv)$ implies
that $R$ corresponds to a {\em regular} (or equitable) partition of $A$, that is, each block of the partition has constant row and column sums. The latter case is of interest to us. 

With regards to the independence number of a graph, let us consider a well-known algebraic bound, namely, the ratio bound (due to Hoffman).
%We now present the inertia bound due to Cvetkovi\'{c}.
%\begin{theorem}(Cvetkovic \label{inertia}\cite{Cvetkovic}) 
%	Suppose the adjacency matrix of a graph $G$ of order $n$ has eigenvalues $\lambda_1\ge \lambda_2\ge \dots\ge \lambda_n$. Assume there are $p$ positive eigenvalues, $q$ negative eigenvalues (counting multiplicities) and the multiplicity of $0$ as an eigenvalue is $z$. Then the size of an independent set is at most $\min\{p + z, q + z\}$. That is,
%	\[\alpha(G)\le |\{i\;:\lambda_i\ge 0\}| \text{ and } \alpha(G)\le |\{i\;:\lambda_i\le 0\}|.\]
%\end{theorem}

%Next is the ratio bound, also called the Hoffman bound.
\begin{theorem}(Hoffman \label{hoffman}\cite{Haemers})
	Let $G$ be a connected  regular graph on $n$ vertices with eigenvalues $\lambda_1\ge \dots \ge \lambda_n$. Then \[\alpha(G)\le n\frac{-\lambda_n}{\lambda_1-\lambda_n}.\] 
\end{theorem}

Abiad, Coutinho, Fiol \cite{abiad2019k} generalized the ratio bound for the independence number to
the $k$-independence number $\alpha_k$. % the following results are known. The first is due to Fiol (see \cite{Fiol}), which involves the so-called $k$-alternating polynomial of $G$. % and will require a preliminary definition. Let $G$ be a graph with distinct eigenvalues $\theta_0> \theta_1>\dots >\theta_d$. Let $P_k(x)$ be chosen among all polynomials $P(x) \in \mathbb{R}_k(x)$, that is, polynomials of real coefficients and degree at most $k$, satisfying $|P(\theta_i)| \le 1$ for all $i =1, ..., d$, and such that $P(\theta_0)$ is maximized. The polynomial $P_k(x)$ defined above is called the $k$-alternating polynomial of $G$ and it was shown to be unique in Fiol, Garriga, and Yebra~\cite{Yebra}, where it was used to study the relationship between the spectrum of a graph and its diameter.
%\begin{theorem}(Fiol\cite{Fiol}) Let $G$ be a $d$-regular graph on $n$ vertices, with distinct eigenvalues $\theta_0> \theta_1>\dots >\theta_d$, and let $P_k(x)$ be its $k$-alternating polynomial. Then,
%	\begin{align}
%	\alpha_k \le \frac{2n}{P_k(\theta_{0})+1}.
%	\end{align}
	%\end{theorem}
%After this, Cvetkovi\'{c}-like and Hoffman-like bounds were shown by Abiad, Cioab\u{a}, and Tait in \cite{Abiad1}. These were further  generalized
 This generalized bound on $\alpha_k$ depends on choosing some polynomial $p(x)
\in \mathbb{R}_k[x]$; that is, $p(x)$ is a polynomial of degree at most $k$ with real coefficients. The proof is based on the interlacing technique described in Theorem~\ref{Interlacing} (see references given there).

%For the rest of this document, let $G$ be a graph on $n$ vertices with adjacency matrix $A$ and eigenvalues $\lambda_1, \lambda_2, \dots, \lambda_n$ such that  $\lambda_1\geq \lambda_2\geq \dots \geq \lambda_n$. We  denote the distinct eigenvalues of $G$ as $\theta_0>\theta_1>\dots > \theta_d$.
Denote $[2,n]=\{2,3,4\dots , n\}$, and let $p(x)\in \mathbb{R}_k[x]$. We define the following parameters:
\begin{itemize}
	\item $W(p)=\max_{u\in V}\{(p(A))_{uu}\}$, where $(p(A))_{uu}$ is the $(u,u)$-entry of the matrix $p(A)$;
	%\item $w(p)=\min _{u\in V}\{(p(A))_{uu}\}$;
	\item $\lambda(p)=\min_{i\in[2,n]}\{p(\lambda_i)\}$.
\end{itemize}

The following is the generalization of the Hoffman bound given by Abiad, Coutinho and Fiol \cite{abiad2019k}.
%\begin{theorem}\cite{abiad2019k}\label{Gbound1} Let $G$ be a graph on $n$ vertices with adjacency matrix $A$ and eigenvalues $\lambda_1, \lambda_2, \dots, \lambda_n$ such that  $\lambda_1\geq \lambda_2\geq \dots \geq \lambda_n$. Let $p(x)\in \mathbb{R}_k[x]$ with corresponding parameters $W(p)$ and $w(p)$. Then the $k$-independence number of $G$ satisfies the bound
%	\begin{align*}
%	\alpha_k\le \min \Big\{|\{i:p(\lambda_i)\ge w(p)\}|,\;|\{i:p(\lambda_i)\le W(p)|\Big\}.
%	\end{align*}
%	
%\end{theorem} 
%
%We now give the Hoffman-like bound.
\begin{theorem}\label{t1}\cite{abiad2019k} 
	Let $G$ be a $\delta$-regular graph on $n$ vertices with adjacency matrix $A$ and eigenvalues $\lambda_1, \lambda_2, \dots, \lambda_n$ such that  $\lambda_1\geq \lambda_2\geq \dots \geq \lambda_n$. Let $p$ be any polynomial in $\mathbb{R}_k[x]$ with corresponding parameters $W(p)$ and $\lambda(p)$, and assume $p(\lambda_1)> \lambda(p)$. Then
	\begin{align}
	\alpha_k\leq n \frac{W(p)-\lambda(p)}{p(\lambda_1)-\lambda(p)}. \label{genb}
	\end{align}
\end{theorem}

\begin{defin}\label{optimal} Let $G$ be a regular graph on $n$ vertices with the eigenvalues $\ev(G)=\{ \theta_0,\theta_1^{m_1},\dots, \theta_d^{m_d} \}$, where $\theta_0>\theta_1>\dots > \theta_d$. 
	Let $B(p,\ev(G)) = n  \frac{W(p)-\lambda(p)}{p(\theta_{0})-\lambda(p)}.$ 
	A polynomial $p\in \mathbb{R}_k[x]$, satisfying $p(\theta_{0})> \lambda(p)$, is called \textit{optimal} with respect to $k$ if it gives the lowest possible bound for $\alpha_k (G)$ of the form $B(p,\ev(G))$.
\end{defin}
Observe that for $c\in \mathbb{R}$, we have $W(p+c)=W(p)+c$, and $W(cp)=cW(p)$. Similarly, $\lambda(p+c)=\lambda(p)+c$ and  $ \lambda(cp)=c\lambda(p)$. Moreover, $(p+c)(\theta_{0})=p(\theta_{0})+c$ and $(cp)(\theta_{0})=cp(\theta_{0})$. Thus,  an optimal polynomial $p(x)\in \mathbb{R}_k[x]$ can be taken to be  monic and without a constant term.
We shall simply write the expression $B(p,\ev(G))$ as $B(p)$ if there is no ambiguity.
As a consequence of Theorem~\ref{t1}, the following results  are known, see Abiad, Coutinho, Fiol \cite{abiad2019k}.
\begin{theorem}\label{t2}\cite{abiad2019k}
	Let $G$ be a $\delta$-regular graph on $n$ vertices  and distinct eigenvalues $\theta_0(=\delta)>\theta_1>\dots > \theta_d$. Let $W_k=W(p)=\max_{u\in V}\{\sum_{i=0}^{k}(A^i)_{uu}\}$. Then, the $k$-independence number of $G$ satisfies the following:
\begin{enumerate}
\item[i.]  If $k=2$, then 
\begin{equation}\label{k2}
\alpha_2\leq n \frac{\theta_0+\theta_i\theta_{i-1}}{(\theta_0-\theta_i)(\theta_0-\theta_{i-1})},
\end{equation}where $\theta_i$ is the largest eigenvalue not greater than $-1$. The optimal polynomial yielding this bound is $p(x)=x^2-(\theta_i+\theta_{i-1})x$.
\item[ii.] If $k>2$ is odd, then %$p(x)=\sum_{j=1}^{k} x^j$ is strictly increasing. This can be seen by observing that the derivative $p'(x)$ can be written as $p'(x)=\sum_{i=1}^{(k-1)/2}i(x^{i-1}+x^i)^2+\frac{k+1}{2}x^{k-1}$, which is a sum of only non-negative terms. So $\lambda(p)=p(\theta_d)$. Hence,  
\begin{align}\label{genk}
\alpha_k\leq n \frac{W_k-\sum_{j=0}^{k}\theta_d^j}{\sum_{j=0}^{k}\delta^j-\sum_{j=0}^{k}\theta_d^j}.
\end{align}
\item[iii.] If $k>2$ is even, then %$p(x)=\sum_{j=1}^{k} x^j$ has  minimum value approaching $-\frac{1}{2}$ as $k\rightarrow \infty$, which implies
\begin{align}
\alpha_k\leq n \frac{W_k+1/2}{\sum_{j=0}^{k}\delta^j+1/2}, 
\end{align}	
%\item [iv.]  If $G =(V, E)$ is a walk-regular graph, then 
% \begin{align}
% \alpha_k\leq n \frac{1-\lambda(q_k)}{q_k(\delta)-\lambda(q_k)}
% \end{align}
% for $k=0,\dots, d-1$, where $q_k=f_0+\dots+f_k$ with the $f_i$'s being the predistance polynomials of $G$ (see \cite{abiad2019k,Fiol,Garriga}),  and $\lambda(q_k)=\min_{i\in [2,d]}\{q_k(\theta_i)\}$. 
where (ii) and (iii) are obtained by using the polynomial $p(x)=\sum_{i=0}^{k}x^i$.
\end{enumerate}
\end{theorem}
 The next result is due to Fiol and it requires a preliminary definition. 
 \begin{defin}
Let $G =(V, E)$ be a graph with adjacency matrix $A$ with the eigenvalues $\ev(G)=\{\theta_0, \theta_1^{m_1},\dots, \theta_d^{m_d}\}$. For a given $k=0,1,\dots, d$, let us consider the set of real polynomials $P_k=\{f\in\mathbb{R}_k[x] :f(\theta_0) =1, f(\theta_i)\ge 0, \text{ for } 1\le i \le d\}$, and the continuous function $\Psi :P_k\rightarrow \mathbb{R}^+$ defined by $\Psi(f)=\tr f(A)$. Then, the \textit{$k$-minor polynomial} of $G$ is the polynomial $f_k$ at which $ \Psi$ attains its minimum:
 	\begin{align}
 	\tr f_k(A) = \min\{\tr f(A) : f \in P_k\}. 
 	\end{align}
 \end{defin}
A graph is called {\em $k$-partially walk-regular} if the number of closed walks of a given length $l\le k$, rooted at a vertex $v$, only depends on $l$. We now state the result by Fiol~\cite{fiol2}.
\begin{theorem}\label{fiolsre}(Fiol \cite{fiol2})
	Let $G$ be a $k$-partially walk-regular graph with n vertices, adjacency matrix A, and the eigenvalues $\ev (G) =\{\theta_0, \theta_1^{m_1},\dots, \theta_d^{m_d}\}$. Let $f_k\in\mathbb{R}_k[x]$ be a $k$-minor polynomial. Then, for every $k=0,\dots , d-1$, the $k$-independence number $\alpha_k$ of $G$ satisfies
	\begin{align}
	\alpha_k\le \tr f_k(A)=\sum_{i=0}^{d} m_if_k(\theta_i).
	\end{align}
\end{theorem}
So, in some sense, a $k$-minor polynomial of $G$ is an optimal polynomial for $\alpha_k(G)$. 
As a consequence of Theorem~\ref{fiolsre}, it was shown by Fiol~\cite{fiol2} that for the case $k=1$, which coincides with the standard independence number, the minor polynomial is $f_i=\frac{x-\theta_d}{\theta_{0}-\theta_d}$. Moreover, $\alpha_1$ coincides with the Hoffman bound in Theorem~\ref{hoffman}. Moreover, for the case $k=2$, the minor polynomial is $f_2=\frac{(x-\theta_i)(x-\theta_{i-1})}{(\theta_{0}-\theta_i)(\theta_{0}-\theta_{i-1})}$, where $\theta_i$  is the largest eigenvalue not greater than $-1$ and $\alpha_2$ is in agreement with Theorem~\ref{t2}$(i)$.
In the general case, $1 \le k\le d$,
 the following was obtained.
 \begin{theorem}\label{fiolsgen}(Fiol~\cite{fiol2})
 	Let $I\subset \{1,\dots, d\}$ range over all index sets with $k$ elements (and if $k$ is odd, it  can be required that $d \in I$). Then,
 	\begin{align}
 	\alpha_k\le \tr f_k(A)=\min_{I}\sum_{j\notin I}m_j\prod_{i\in I}\frac{\theta_j-\theta_i}{\theta_0-\theta_i}.
 	\end{align}
 \end{theorem}
  
For the case $k=3$ in particular, Fiol proposed the minor polynomial $f_3=f_1f_2$ as good (and often optimal) choice. With this choice of polynomial we have the following.
\begin{corollary}\label{fiolk33}(Fiol \cite{fiol2}) If $G$ is at least $3$-partially walk-regular, and  $n_t$ is the common number of triangles rooted at every vertex of $G$, then 
	\begin{align}\label{fiolk333}
	\alpha_3\le \tr f_3 = n\frac{2n_t-\theta_{0}(\theta_d+\theta_i+\theta_{i-1})-\theta_d\theta_i\theta_{i-1}}{(\theta_0-\theta_d)(\theta_{0}-\theta_i)(\theta_{0}-\theta_{i-1})},
	\end{align} where $\theta_i$  is the largest eigenvalue not greater than $-1$.
\end{corollary}

\section{The case $k$=3 (the $3$-independence number)}
%\subsection{The Polynomial of the form $p(x)=x^3+bx^2+cx$ }
 
In the next theorem, we obtain an optimal polynomial for $k=3$, and hence an optimal bound for $\alpha_3$ over all bounds of the form given in  Theorem~\ref{t1}. The following Lemma will be needed. %For the purpose of generality, we include $p(\theta_0)$ in the definition of $\lambda(p)$ just in this lemma.

If we consider $p(x)=x^3+bx^2+cx$ and set $x=\theta_i$, and let $c$ vary, we get a linear function $p_i(c)=\theta_i^3+b\theta_i^2+c\theta_i$ of slope $\theta_i$. Since $\theta_d< \theta_i$ for all $i\neq d$ then for $c$ sufficiently large enough, $p_d(c)<p_i(c)$ and for $c$ small $p_d(c)>p_i(c)$. 
%Now let ${c^*_i}$ be the $c$ coordinate of intersection of $p_i(c)$ and $p_{i-1}(c)$ and define  $\Upsilon(c) =\min \{p_i(c): i\ne d\}$. Then  $\Upsilon(c)$ is a concave down piece-wise linear function,	and $\Upsilon(c)=p_i(c)$  if $c\in [c_{i}^*,c_{i+1}^*].$ Observe that, at $c=\frac{b^2}{3}$, we have $p_{d}(c)< p_{m}(c)<p_{{m-1}}(c)<\dots<p_{1}(c)$ for some $m\le d$. Hence,  the line $p_d(c)$ has the least slope and  will intersect $\Upsilon(c)$ at some point say, $(c^*, \Upsilon(c^*))$. Then, $c^*$ is the largest value  of $c$ such that $p_d(c)=p_i(c)$ for some value of $i$. Thus,  $p_d(c)\le \Upsilon(c)$  for $c\in [c^*, \frac{b^2}{3})$.

\begin{lemma}\label{lemma5}
	Let  $d\geq3$, and $\theta_0, \theta_1,\dots ,\theta_d$ be real numbers with $\theta_{0}>\dots>\theta_d$.
	Consider $b\in \mathbb{R}$ fixed throughout. Let  $p(x)=x^3+bx^2+cx$, where $c\in \mathbb{R}$. For $i=0,1, \dots,d$, let  $p_i(c)=p(\theta_i)$ and $\lambda(p)=\min\{ p_i(c):i=1,\dots,d\}$. 
	 Let $c^*$ be the largest value  of $c$ such that $p_d(c)=p_i(c)$ for some value of $i$, and denote such $i$ as $j$.	
	Let $c=c_i^*$ be the solution of the linear equation $p_i(c)=p_{i-1}(c)$ for $i=1,2,\dots, j$.
	%the largest $j\in \{1,\dots,m\}$ such that ${c_{j}^*}<  {c^*}$
	  Define the intervals
	\begin{align*}%I_0&=(-\infty,c_0^*]
	 I_s&=[c_{s}^*,c_{s+1}^*] \text{ for } s=1,\dots ,j-1, \quad
	 I_j= [c^*_j,c^*], \text{ and } 
	I_d =[c^*,\infty).
	\end{align*}
	Then, $\lambda(p)=p(\theta_i)$ if $c\in I_i$ for $i\in \{1,\dots, j\}\cup \{d\}$.%, that is, $\theta_i$ minimizes the polynomial $p$ over the set of eigenvalues $\theta_1,\dots,\theta_j ,\theta_d$ if  $c\in I_i$. 	
\end{lemma}
\begin{proof}Let $\Lambda=\{ \theta_{0}, \dots, \theta_{d} \}$, so $\lambda(p)=\min\{ p(x):x\in \Lambda\backslash(\theta_{0})\}$.
	
	If we differentiate $p(x)$, we observe that it has a local maximum at $x_1=\frac{1}{3}(-b-\sqrt{b^2-3c})$ and a local minimum at $x_2=\frac{1}{3}(-b+\sqrt{b^2-3c})$ if $c<\frac{b^2}{3}$. If $c\ge \frac{b^2}{3}$,  then $p(x)$ is increasing and has no turning point and  hence, on $[\theta_d, \infty ]$ it attains an absolute minimum at $ \theta_d$. Thus, we have $\lambda(p)=p(\theta_d) \text{ if } c\ge\frac{b^2}{3}$.

	So, for $c<\frac{b^2}{3}$, we have that on $[\theta_d, \infty )$, $p(x)$ attains an absolute minimum at $x_2$ or at $\theta_d$. We explore this case further.
	
	Consider $x_2$ as a function of $c$; that is, $x_2(c)=\frac{1}{3}(-b+\sqrt{b^2-3c})$. Observe that $x_2(c)$ is  continuous and decreasing on $(-\infty,\frac{b^2}{3})$, with $\lim_{c\rightarrow \frac{b^2}{3}}x_2(c)=\frac{-b}{3}$ and $\lim_{c\rightarrow -\infty}x_2(c)=\infty$. 
	So, as $c$ increases, the local minimum occurs at smaller values of $x_2$.
	Let $m$ be  the maximum index such that $\theta_m> \frac{-b}{3}$. Then, over $\Lambda$, for $i=0,1,\dots,m$, there exists $c_i\in (-\infty,\frac{b^2}{3})$ such that $x_2(c_i)=\theta_i$, and $c_0<c_1<\dots<c_m$. That is, if $c=c_i$, then $x_2=\theta_i$ and absolute minimum of $p(x)$ occurs at $x=\theta_i$ or $x=\theta_d$.
	
	%For $i\in \{1,\dots,m\}$, we recall ${c^*_i}$ is the $c$ coordinate of the intersection of $p_i(c)$ and $p_{i-1}(c)$. Since $c_0<c_1<\dots < c_m$ and $c_{i-1}< {c_i^*}<  c_{i}$, we have ${c_1^*}< {c_2^*}< \dots < {c_{m}^*} < {c_m}$. 	Now since $\theta_{i}<\theta_{i-1}$,  if $c \ge  {c_{i}^*}, \text{ then }p_{i}(c)<p_{i-1}(c)<\dots < p_1(c) $. Similarly, if   $c \le  {c_{i+1}^*}, \text{ then }p_{i}(c)<p_{i+1}(c)<\dots < p_m(c)$. Hence, if $c\in [c_{i}^*,c_{i+1}^*]$, then $\Upsilon(c)=\min\{ p_s(c):s\in \{1,2,\dots, m\}\}=p_{i}(c)$. Now line $p_d(c)$, having the least slope,  intersects $\Upsilon(c)$ at the point $(c^*, \Upsilon(c^*))$. %Then $c^*$ is the largest value  of $c$ such that $p_d(c)=p_i(c)$ for some value of $i$. 	And so  $p_d(c)\le \Upsilon(c)$  for $c\in [c^*, \frac{b^2}{3})$.
	Now, for $i\in \{1,\dots,m\}$, let ${c^*_i}$ be the $c$ coordinate of intersection of $p_i(c)$ and $p_{i-1}(c)$ and define  $\Upsilon(c) =\min \{p_i(c)\}$. Then  $\Upsilon(c)$ is a concave down piece-wise linear function. Since $c_0<c_1<\dots < c_m$ and $c_{i-1}< {c_i^*}<  c_{i}$, we have ${c_1^*}< {c_2^*}< \dots < {c_{m}^*} < {c_m}$. 
	Now since $\theta_{i}<\theta_{i-1}$,  if $c \ge  {c_{i}^*} \text{ then }p_{i}(c)<p_{i-1}(c)<\dots < p_1(c) $. Similarly, if   $c \le  {c_{i+1}^*} \text{ then }p_{i}(c)<p_{i+1}(c)<\dots < p_m(c)$. Hence,  $\Upsilon(c)=p_i(c)$  if $c\in [c_{i}^*,c_{i+1}^*].$ Observe that, at $c=\frac{b^2}{3}$, we have $p_{d}(c)< p_{m}(c)<p_{{m-1}}(c)<\dots<p_{1}(c)$. Hence, line $p_d(c)$, having the least slope, intersects $\Upsilon(c)$ at some point say, $(c^*, \Upsilon(c^*))$. Then, $c^*$ is the largest value  of $c$ such that $p_d(c)=p_i(c)$ for some value of $i$. Thus,  $p_d(c)\le \Upsilon(c)$  for $c\in [c^*, \frac{b^2}{3})$
	
	Hence, for the largest $j\in \{1,\dots,m\}$ such that ${c_{j}^*}<  {c^*}$, we have ${c_1^*}< {c_2^*}< \dots < {c_{j}^*} < {c^*}$, and so %$\lambda(p)=p_s(c) \text{ if } c\in [c_{s}^*,c_{s+1}^*]$ for $s\in \{1,\dots,j-1\}$, $\lambda(p)=p_j(c) \text{ if } c\in [c_{j}^*,c^*]$ and $\lambda(p)=p_d(c) \text{ if } c\in [c^*,\infty)$.
	\[\lambda(p)=
	\begin{cases}
		p_s(c), \text{ if }  c_{s}^*\le c\le c_{s+1}^* \text{ for all } 1\le s\le j-1,\\
			p_j(c), \text{ if } c_{j}^*\le c \le c^*,\\
			p_d(c), \text{ if }  c^*\le c < \infty.
	\end{cases}\]
\end{proof}

We now present our main result.
\begin{theorem}\label{Gmaink3} Let $G$ be a $\delta$-regular graph with $n$ vertices, adjacency matrix $A$, and distinct eigenvalues $\delta=\theta_0>\theta_1>\dots> \theta_d$, with $d\geq3$. 	Let $s$ be the largest index such that $\theta_s\ge -\frac{\theta_{0}^2+\theta_{0}\theta_d-\triangle}{\theta_{0}(\theta_d+1)}$, where $\triangle=\text{max}_{u\in V}\{(A^3)_{uu}\}$. %, twice the largest number of triangles on any vertex in $G$. 
	Let $b=-(\theta_s+\theta_{s+1}+\theta_d)$
	and  $c=\theta_d\theta_s+\theta_d\theta_{s+1}+\theta_s\theta_{s+1}$.  
	Then, $p(x)=x^3+bx^2+cx$ is an optimal polynomial for $k=3$. The corresponding bound on the $3$-independence number of $G$ is
	\begin{align}\label{Gmaink-3}
	\alpha_3%&\le n\frac{\triangle-\theta_{d}^3 +b(\theta_{0}-\theta_d^2)-c\theta_{d} }{(\theta_{0}^3-\theta_d^3)+b(\theta_{0}^2-\theta_d^2)+c(\theta_{0}-\theta_{d}) }\nonumber\\
	&\le n\frac{\triangle-\theta_{0}(\theta_s+\theta_{s+1}+\theta_d)-\theta_{s}\theta_{s+1}\theta_{d} }{(\theta_{0}-\theta_s)(\theta_{0}-\theta_{s+1})(\theta_{0}-\theta_{d})}.
	\end{align} 
	If equality is attained, then the matrix $A^3-(\theta_s+\theta_{s+1}+\theta_d)A^2+(\theta_d\theta_s+\theta_d\theta_{s+1}+\theta_s\theta_{s+1})A$ has a regular partition (with a set of $\alpha_3$ $3$-independent vertices and its complement) with quotient matrix
	\begin{align}\label{B3mat}
	B=\left[\begin{array}{cc}
	B_{11} & B_{12}\\ B_{21}& B_{22}
	\end{array}\right],
	\end{align}
	where % \begin{gather*}
	{\footnotesize \begin{align*}
	B_{11}&= \triangle-(\theta_s+\theta_{s+1}+\theta_d)\delta,\\
	B_{12}&=\delta^3-(\theta_s+\theta_{s+1}+\theta_d)\delta^2 + (\theta_d\theta_s+\theta_d\theta_{s+1}+\theta_s\theta_{s+1}+\theta_s+\theta_{s+1}+\theta_d)\delta-\triangle,\\
	B_{21}&= \triangle-(\theta_s+\theta_{s+1}+\theta_d)\delta-\theta_s\theta_{s+1}\theta_d,\\
	B_{22}&= \delta^3-(\theta_s+\theta_{s+1}+\theta_d)\delta^2 + (\theta_d\theta_s+\theta_d\theta_{s+1}+\theta_s\theta_{s+1}+\theta_s+\theta_{s+1}+\theta_d)\delta+\theta_s\theta_{s+1}\theta_d-\triangle.
	\end{align*}}
	%\end{gather*} 
\end{theorem}

If we assume $G$ is at least $3$-partially walk-regular, and let $n_t$ be the common number of triangles rooted at every vertex of $G$, then $\triangle=2n_t$ and so
\begin{align}\label{Gmaink-33}
\alpha_3 &\le n\frac{2n_t-\theta_{0}(\theta_s+\theta_{s+1}+\theta_d)-\theta_{s}\theta_{s+1}\theta_{d} }{(\theta_{0}-\theta_s)(\theta_{0}-\theta_{s+1})(\theta_{0}-\theta_{d})},
\end{align} 
which is consistent with Fiol's bound in [\ref{fiolk33}] if its index $i$ coincides with index $s+1$, otherwise our bound is stronger. We will revisit this in Section~\ref{fiolcompare}.

In particular, if $G$ is a bipartite graph, then since $\theta_d=-\theta_0$ and $\triangle=0$, we have that $\theta_s$ is the least eigenvalue greater or equal to $0$, and the corresponding bound on $\alpha_3$ in (\ref{Gmaink-3})  becomes \begin{align}\label{Gmaink-3co}
\alpha_3%&\le n\frac{(\theta_{0}^2+c) +b(1-\theta_0) }{2(\theta_{0}^2+c) }\nonumber\\
&\le \frac{\theta_{0}-(\theta_{s}+\theta_{s+1})+\theta_{s}\theta_{s+1}}{2(\theta_{0}-\theta_{s})(\theta_{0}-\theta_{s+1})}.
\end{align}
We now prove Theorem~\ref{Gmaink3}.

\begin{proof}  
	Take any $b,c\in \mathbb{R}$, and define a polynomial $p(x)=x^3+bx^2+cx$. We show that  $b$ and $c$ as stated in the theorem give an optimal polynomial.  Assume  $\lambda(p)$ occurs at $\theta_{t}$.  
	Then, %$\lambda(p)=p(\theta_{t})=\theta_{t}^3+b\theta_{t}^2+c\theta_{t}.$  We have $p(\theta_0)=\theta_{0}^3+b\theta_{0}^2+c\theta_{0}$, and $W(p)=\text{max}_{u\in V}\{(p(A))_{uu}\} =\text{max}_{u\in V}\{(A^3)_{uu}+b\,(A^2)_{uu} +c\,(A)_{uu}\}=\triangle+b\,\delta+c(0)=\triangle+b\,\theta_{0}$. 
	$W(p)=\triangle+b\,\theta_{0}$.
	Thus, the bound in Equation~(\ref{genb}) becomes
	\begin{align}\label{keyg}
	\varPhi_t(b,c)&=n\frac{\triangle+b\theta_{0}-p(\theta_t)}{p(\theta_{0})-p(\theta_t)}= n\frac{\triangle+b\theta_{0}-(\theta_{t}^3+b\theta_{t}^2+c\theta_{t})}{\theta_{0}^3+b\theta_{0}^2+c\theta_{0}-(\theta_{t}^3+b\theta_{t}^2+c\theta_{t})}.
	\end{align} 
	We investigate the pairs $(b,c)$ that minimize the  bound $\varPhi_t(b,c)$. %To satisfy the condition $\lambda(p_{b,c})=p_{b,c}(\theta_i)< p_{b,c}(\theta_0)$ in Definition~\ref{optimal}, we must have that $b(\theta_i+\theta_0)+c > -(\theta_i^2+\theta_i\theta_0+\theta_0^2)$.
	Let $b\in \mathbb{R}$ be fixed. Differentiating (\ref{keyg}) with respect to $c$ gives
	\begin{align*}
	\frac{\partial}{\partial c}\varPhi_t(b,c)&= n\frac{ \theta_{0} \theta_{t}^{3}+b \theta_{0} \theta_{t}^{2}  + {\left(-b \theta_{0}^{2} - \theta_{0}^{3} +b \theta_{0} + \triangle\right)} \theta_{t} -b \theta_{0}^{2} - \theta_{0} \triangle}{\big[\theta_{0}^3+b\theta_{0}^2+c\theta_{0}-(\theta_{t}^3+b\theta_{t}^2+c\theta_{t})\big]^2}\\
	&= n\frac{(\theta_{0}-\theta_{t})(-b \theta_{0} \theta_{t} - \theta_{0}^{2}\theta_t -  \theta_{0}\theta_{t}^{2} -b \theta_{0} - \triangle)}{ (\theta_{0}-\theta_{t})^2(b \theta_{0} + \theta_{0}^{2} +b\theta_{t} + \theta_0 \theta_{t} +\theta_{t}^2+c)^2}\\
	&= -n\frac{\theta_{0}\theta_{t}^{2}+(b \theta_{0}  + \theta_{0}^{2})\theta_t   + (b\theta_{0} + \triangle)}{ (\theta_{0}-\theta_{t})(b \theta_{0} + \theta_{0}^{2} +b\theta_{t} + \theta_0 \theta_{t} +\theta_{t}^2+c)^2}.
	\end{align*} 
	The denominator is equal to $(p(\theta_{0})-p(\theta_t))^2$ divided by $\theta_0-\theta_t$. Since $p(\theta_{0})-p(\theta_t)>0$ and $\theta_0-\theta_t>0$, the denominator is positive. The numerator is a quadratic function in $x=\theta_t$,  that is,  $f(x)=\theta_{0}x^{2}+(b \theta_{0}  + \theta_{0}^{2})x   + (b\theta_{0} + \triangle)$ with roots 
	
	%$$\left[\theta_{i} =\frac{-b \theta_{0} - \theta_{0}^{2} \pm \sqrt{(-b \theta_{0} - \theta_{0}^{2})^2+4 \,\theta_{0}(-b \theta_{0} - \triangle) }}{2 \, \theta_{0}}\right]$$
	$$x_{1,2} =-\frac{b+\theta_{0}}{2}\pm \sqrt{\left(\frac{b+\theta_{0}}{2} \right)^2-b- \frac{ \triangle}{\theta_{0}}}.$$ 
	Let $\nu=-\frac{b+\theta_{0}}{2}- \sqrt{\left(\frac{b+\theta_{0}}{2} \right)^2-b- \frac{ \triangle}{\theta_{0}}}$ and $\mu=-\frac{b+\theta_{0}}{2}+ \sqrt{\left(\frac{b+\theta_{0}}{2} \right)^2-b- \frac{ \triangle}{\theta_{0}}}$. 
	Since each vertex of $G$ can be on at most $ \binom{\theta_0}{2}$ different triangles, we have $\triangle\le \theta_0(\theta_0-1)$.	So, we have \begin{align*}
	\left(\frac{b+\theta_{0}}{2} \right)^2-b- \frac{ \triangle}{\theta_{0}}&\ge \left(\frac{b+\theta_{0}}{2} \right)^2-b- \frac{ \theta_{0}(\theta_{0}-1)}{\theta_{0}} \\
	&= \left(\frac{b+\theta_{0}}{2} \right)^2-b- (\theta_{0}-1)= \left(\frac{b+\theta_{0}}{2} \right)^2-(b+ \theta_{0})+1\\
	&= \left(\frac{b+\theta_{0}}{2}-1\right)^2\ge 0. 
	\end{align*}
	Thus, $\nu$ and $\mu$ are real numbers, and
	 \begin{itemize}
		\item $\varPhi_t(b,c)$ is increasing with respect to $c$ if 
		$\nu<  \theta_{t} <\mu$,
		\item $\varPhi_t(b,c)$ is decreasing with respect to $c$ if  $ \theta_{t} < \nu$ or $\theta_{t} >\mu$, 
		\item $\varPhi_t(b,c)$ is constant with respect to $c$ if  $\theta_{t} =\nu \text{ or } \theta_t=\mu$.
	\end{itemize} 
	
	% Since $c$ plays no direct role (after $\theta_i$ is fixed) in where the bound is increasing or decreasing, we can fix $b$ and then determine an optimal $c$. 

	%Solving $\theta_i\le \nu$ or $\theta_i\ge \mu$ we have that  $b\ge \frac{-1}{1+\theta_{i}}(\theta_i\theta_{0}+\frac{\triangle}{\theta_{0}}+\theta_{i}^2)$ if $\theta_i >-1 $ and $b\le \frac{-1}{1+\theta_{i}}(\theta_i\theta_{0}+\frac{\triangle}{\theta_{0}}+\theta_{i}^2)$ if $\theta_i<-1$. Similarly, if $\nu\le \theta_i \le \mu$, then $b\le \frac{-1}{1+\theta_{i}}(\theta_i\theta_{0}+\frac{\triangle}{\theta_{0}}+\theta_{i}^2)$ if $\theta_i >-1 $ and $b\ge \frac{-1}{1+\theta_{i}}(\theta_i\theta_{0}+\frac{\triangle}{\theta_{0}}+\theta_{i}^2)$ if $\theta_i<-1$. Using the fact that $\nu \le -1$ and $\mu\ge -1$, we have the following cases.
 Now we show $\nu\le -1\le \mu$.
	\begin{align*}\mu&=-\frac{b+\theta_{0}}{2}+ \sqrt{\left(\frac{b+\theta_{0}}{2} \right)^2-b- \frac{ \triangle}{\theta_{0}}}\\&\ge -\frac{b+\theta_{0}}{2}+\sqrt{\left(\frac{b+\theta_{0}}{2}-1\right)^2}\\
	&= -\frac{b+\theta_{0}}{2}+\left|\frac{b+\theta_{0}}{2}-1\right|\\
	&= \begin{cases}
	-\frac{b+\theta_{0}}{2}+(-(\frac{b+\theta_{0}}{2}-1)) & \text{ if } b\le 2-\theta_{0}\\
	-\frac{b+\theta_{0}}{2}+(\frac{b+\theta_{0}}{2}-1) & \text{ if } b\ge 2-\theta_{0}
	\end{cases}\\
	&= \begin{cases}
	-(b+\theta_{0})+1\ge-1 & \text{ if } b\le 2-\theta_{0}\\
	-1 & \text{ if } b\ge 2-\theta_{0}
	\end{cases}.
	\end{align*}
	Hence, we have $\mu\ge -1$.
	
	Also, we have 
	\begin{align*}
	\nu &= -\frac{b+\theta_{0}}{2}- \sqrt{\left(\frac{b+\theta_{0}}{2} \right)^2-b- \frac{ \triangle}{\theta_{0}}}\\
	&\le -\frac{b+\theta_{0}}{2}-\sqrt{\left(\frac{b+\theta_{0}}{2}-1\right)^2}\\
	&= -\frac{b+\theta_{0}}{2}-\left|\frac{b+\theta_{0}}{2}-1\right|\\
	&= \begin{cases}
	-\frac{b+\theta_{0}}{2}-(-(\frac{b+\theta_{0}}{2}-1)) & \text{ if } b\le 2-\theta_{0}\\
	-\frac{b+\theta_{0}}{2}-(\frac{b+\theta_{0}}{2}-1) & \text{ if } b\ge 2-\theta_{0}
	\end{cases}\\
	&= \begin{cases}
	-1 & \text{ if } b\le 2-\theta_{0}\\
	-(b+\theta_{0})+1\le -1 & \text{ if } b\ge 2-\theta_{0}
	\end{cases}	.
	\end{align*} 
	Thus, we have $\nu\le -1$.
	
	Now, for a fixed $b$, we determine an optimal $c$, that is, a $c$ that will minimize $\varPhi_t(b,c)$. 	
Adopt the notation of Lemma~\ref{lemma5}, and  let  $p_i(c)=p(\theta_i)$ and  $\lambda(p)=\min\{ p_i(c):i=1,\dots,d\}$, we define the intervals \[I_\ell=[c_{\ell},c_{\ell+1}] \text{ for } \ell=1,\dots ,j-1,\;I_j=[c_{j},c^*],  \text{ and }   I_d=[c^*,\infty).\] Recall that $c^*$ is the largest value of $c$ such that $p_d(c)=p_i(c)$ for some value of $i$ denoted as $j$, and $c=c_i$ is the solution of the linear equation $p_i(c)=p_{i-1}(c)$ for $i=1,\dots, j$. Recall also that by Lemma~\ref{lemma5}, if $c\in I_i$, then $\lambda(p)=p(\theta_i)$, that is, $\theta_i$ minimizes $p$. Observe that, by choosing any $c\in I_i$ for $i\in \{1,\dots,j\}\cup \{d\}$, we ensured $p(\theta_0)>\lambda(p)$.  Note that, for $b$ fixed, we aim to find $\Omega=\min\{\varPhi_i(b,c):i\in \{1,2,\dots, j,d\}, c\in I_i \}$. As $\lambda(p)$ is constant and $\varPhi_i(b,c)$ is monotone on each of the intervals $I_1,\dots ,I_j,I_d$, we only need to compare the values of $\varPhi_i(b,c)$ at interval endpoints. In particular, since $\nu\le -1\le \mu$, we observe the following.
Let  $i'$ be such that  $\theta_{i'}$ is the largest eigenvalue with $\theta_{i'}\le \mu$, and let $j'$ be such that $\theta_{j'}$ is the largest eigenvalue with $\theta_{j'}\le \nu$.
	\begin{itemize}
		\item[Case 1.] Assume $\theta_i\le \nu$. In this case, $\varPhi_i(b,c)$ is decreasing or constant on $I_i$. 
		Hence, if $i\ne d$, then the minimum of $\varPhi_i(b,c)$ exists and is attained at the right endpoint of the interval, that is, at 
		$$\begin{cases}
		c_{i+1} & \text{ if } j' \le i \le j-1,\\
		c^* & \text{ if } i=j.
		\end{cases}$$ In particular,
		we have
		\[ \varPhi_m(b,c_{m+1})=\varPhi_{m+1}(b,c_{m+1})\ge\varPhi_{m+1}(b,c_{m+2}) \text{ for all } j'\le m\le j-1, \] and 
		\[\varPhi_{j-1}(b,c_{j})=\varPhi_{j}(b,c_{j})\ge \varPhi_j(b,c^*)=\varPhi_{d}(b,c^*)\ge\varPhi_{d}(b,c)\text{ for all }  c > c^*.\]
		
		Thus, we have \[\Omega=\min\{\varPhi_i(b,c):i\in \{j', j'+1,\dots,j,d\}, c\in I_i \}=\varPhi_{d}(b,c) \text{ for all }  c > c^*.\]
		However, $\varPhi_{d}(b,c)$ has no minimum on $I_d$ but
		$\lim\limits_{c\rightarrow \infty}\varPhi_d(b,c)=n\frac{-\theta_d}{\theta_{0}-\theta_{d}}$ is an upper bound on $\alpha_3 $ and it is less than any other bound in this case. But this is the Hoffman bound for $\alpha_1$, which clearly holds for $\alpha_3$.
		
		\item[Case 2.] Assume  $\mu<  \theta_i$. In this case, $\varPhi_i(b,c)$ is decreasing on $I_i$. Hence, $\varPhi_i(b,c)$ attains its minimum at the right endpoint of $I_i$; that is, at $c_{i+1}$ for $1< i \le i'-1.$ Recall that $\theta_{i'-1}$ is the smallest eigenvalue greater than $\mu$. In particular, the minimum occurs at		
		$$\begin{cases}
		c_{i+1} & \text{ if } 1 \le i \le j-1,\\
		c^* & \text{ if } i=j.
		\end{cases}$$
		
		We observe the following. If $m\in \{1, \dots i'-1 \}$, then
		\begin{align*}
		\varPhi_m(b,c_{m+1})
		&= \varPhi_{m+1}(b,c_{m+1})> 
		\varPhi_{m+1}(b,c_{m+2}) \text{ for all } 1\le m \le i'-1.
		\end{align*}	
		
		Thus, we have \[\Omega=\min\{\varPhi_i(b,c):i\in \{1, \dots,i'-1\}, c\in I_i \}=\varPhi_{i'-1}({b,c_{i'}}).\]	
		
		Note that $j>i'-1$ or $j\le i'-1$. %Thus, $c_{\jmath-1}=c_{k}$	if  $\jmath-1=k$. 
		Thus, we have two possibilities here, that is,  either $\theta_j< \mu$ or $\theta_j\ge \mu$.
		In particular,
		$$\Omega = \begin{cases}
		\varPhi_{i'-1}(b,c_{i'}) & \text{ if } j\ge  i',\\
		\varPhi_{j}(b,c^*) & \text{ if }  j< i'.
		\end{cases}$$
		That is, for all $i$ such that $\mu<\theta_i$, the bound $\varPhi_i(b,c)$ is minimized at $i=i'-1$ with $c=c_{i'}$ if $j\ge i'$, or at $i=j$ with $c=c^*$ if $j\le i'-1$.
		\item[Case 3.] Assume $\nu\le \theta_i\le \mu$.
		In this case, $\varPhi_i(b,c)$ is increasing or constant on $I_i$.
		Hence, $\varPhi_i(b,c)$ attains its minimum at the left endpoint of $I_i$; that is, at
		$$\begin{cases}
		c_{i} & \text{ if } i'\le i \le j,\\
		c^* & \text{ if } i=d.
		\end{cases}$$
		
		Now we have that  $\varPhi_{d}(b,c^*)= \varPhi_{j}(b,c^*)\ge\varPhi_{j}(b,c_{j})$. Moreover,
		\begin{align*}
		\varPhi_m(b,c_{m})= \varPhi_{m-1}(b,c_{m})\ge 
		\varPhi_{m-1}(b,c_{m-1}) \text{ for all } m\le j.
		\end{align*}
		Thus, we have \[\Omega=\min\{\varPhi_i(b,c):i\in \{i',\dots,j,d \}, c\in I_i \}=\varPhi_{i'}(b,c_{i'}).\]
		
		Observe that we have two possibilities here also, that is,  either $\theta_j< \mu$ or $\theta_j\ge \mu$. In particular,
		$$\Omega = \begin{cases}
		\varPhi_{i'}(b,c_{i'}) & \text{ if } j\ge i',\\
		\varPhi_{d}(b,c^*) & \text{ if } j< i'.%,\text{ i.e., } k=i'-1.
		\end{cases}$$
		That is, for all $i$ such that $\nu\le \theta_i\le \mu$, the bound $\varPhi_i(b,c)$ is minimized at $i=i'$ with $c=c_{i'}$ if $j\ge i'$ or at $i=d$ with $c=c_*$ if $j\le i'-1$.
	\end{itemize}
 
	Note that $\varPhi_{i'}({b,c_{i'}})=\varPhi_{i'-1}({b,c_{i'}})$ and $\varPhi_{j}(b,c^*)=\varPhi_{d}(b,c^*)$. Thus, case 2 and case 3 yield the same bounds. That is, if $j\ge i'$, we have $$\Omega =\varPhi_{i'}({b,c_{i'}})=\varPhi_{i'-1}({b,c_{i'}})$$ and if $j< i'$, we have $$\Omega= \varPhi_{j}(b,c^*)=\varPhi_{d}(b,c^*).$$ 	 
	    %Moreover, $c=c_{\jmath-1}^*$ is  the solution of the linear equation $p_{b,c}(\theta_{\jmath-1})=p_{b,c}(\theta_{\jmath})$.
Now, let's investigate the two bounds arising from Cases 2 and 3.
	\begin{itemize}
		\item [Case A:]  $j\ge i'$, that is,  $\theta_j\le \mu $. Then,
		$$\Omega=\varPhi_{i'}({b,c_{i'-1}})=\varPhi_{i'-1}({b,c_{i'-1}}) .$$
		
		Note that $\mu$ is not necessarily an eigenvalue of $G$. But we have the bound  $\Omega=\varPhi_{i'}({b,c_{i'-1}})=\varPhi_{i'-1}({b,c_{i'-1}}) \ge \lim\limits_{c\rightarrow \infty}\varPhi(\mu,b,c)=n\frac{-\mu}{\theta_{0}-\mu}$, and so $-1\le \mu<0$. We have equality when $\theta_{i'}=\mu$, and hence, $n\frac{-\mu}{\theta_{0}-\mu}$ becomes a bound on $\alpha_3$.
		
		% [ or use the fact that $\omega_{\jmath}(b,c)$ is minimized at approximately the constant value where $\theta_{\jmath}=\mu$. That is, $\omega_{\jmath-1}({b,c_{\jmath-1}^*})\ge \omega_{\ell}({b,c_{}})$ where $\mu=\theta_{\ell}.$ Hence $\min_{c\in \RR} \omega_i(b,c)\ge \lim\limits_{c\rightarrow \infty}\omega_\ell(b,c)=n\frac{-\theta_\ell}{\theta_{0}-\theta_{\ell}}$ and $-1\le \theta_\ell<0$. This is not completely accurate either. Need more work; to be revised.]
		
		\item  [Case B:] $j< i'$, that is,  $\mu \le \theta_j< \theta_{0}$. 		We have  \[\Omega=\varPhi_{d}(b,c^*)=\varPhi_{j}(b,c^*).\]
		We will investigate this  further.
	
	\end{itemize}
	Local minimum of $p(x)$ occurs at  $x=\tau\ge \mu$. So, if
	$$\lambda(p)=p(\theta_{d})=p(\tau),$$
	then \begin{align}c=-(\tau^2+\theta_d\tau+\theta_d^2)-b(\theta_d+\tau).\label{c}
	\end{align}  
	Substituting (\ref{c}) into Equation~(\ref{keyg}), we have  %$$\lambda(p)=p(\tau)=\tau^3+b\tau^2+c\tau=-(b\tau\theta_d+\tau\theta_d^2+\tau^2\theta_d).$$ So Equation~(\ref{keyg}) becomes 
	\begin{align*}
	\varPhi(b)%&= n\frac{\triangle+b\theta_{0}+(b\tau\theta_d+\tau\theta_d^2+\tau^2\theta_d) }{\theta_{0}^3+b\theta_{0}^2-\theta_{0}(\tau^2+\theta_d\tau+\theta_d^2)-b\theta_{0}(\theta_d+\tau)+(b\tau\theta_d+\tau\theta_d^2+\tau^2\theta_d)}\\
	&= n\frac{\triangle+b\theta_{0}+(b\tau\theta_d+\tau\theta_d^2+\tau^2\theta_d) }{(\theta_{0}+b+\tau+\theta_d)(\theta_{0}-\tau)(\theta_{0}-\theta_d)} ,
	\end{align*} 
	and differentiating with respect to $b$ gives
	\begin{align*}
	\varPhi'(b)%&= n\frac{(\theta_{0}\theta_d\tau+\theta_{0}\tau+\theta_{0}^2+\theta_{0}\theta_d-\triangle)(\theta_{0}-\tau)(\theta_{0}-\theta_d)}{\big[(\theta_{0}+b+\tau+\theta_d)(\theta_{0}-\tau)(\theta_{0}-\theta_d) \big]^2}\\
	&= n\frac{(\theta_{0}\theta_d+\theta_{0})\tau+(\theta_{0}^2+\theta_{0}\theta_d-\triangle)}{(\theta_{0}+b+\tau+\theta_d)^2(\theta_{0}-\tau)(\theta_{0}-\theta_d) }.
	\end{align*} 
	We are interested in where $\varPhi'(b)=0$, so  $\tau=-\frac{\theta_{0}^2+\theta_{0}\theta_d-\triangle}{\theta_{0}(\theta_d+1)}$. 
	
	Let $\theta_s$ be an eigenvalue such that $\lambda(p)=p(\theta_s)$. Then, $\theta_s$ is the least eigenvalue of $G$ such that $\theta_s\ge \tau$ or the largest eigenvalue of $G$ such that $\theta_s\le \tau$. Hence, we want $b$ and $c$ that will give $\lambda(p)=p(\theta_d)=p(\theta_s)$, this being a necessary condition. So, we have \begin{align}c=-(\theta_s^2+\theta_d\theta_s+\theta_d^2)-b(\theta_d+\theta_s).\label{cc}
	\end{align}  We now proceed to find $b$.
	
	We note that, if $\tau=\theta_s=-\frac{\theta_{0}^2+\theta_{0}\theta_d-\triangle}{\theta_{0}(\theta_d+1)}$, then $\varPhi(b)$ is a constant function of $b$. Thus any value of $b$ that ensures the minimum still occurs at $\theta_s$ will give us an optimal bound. That is, the domain for $b$ must be such that, $\lambda(p_{b})$ will always occur at $\theta_s$. Now let us solve for the domain of $b$. We have that $b$ satisfies $p_{b}(\theta_s)\le p_{b}(\theta_{s-1})$ with strict inequality if $s=1$, and $p_{b}(\theta_s)\le p_{b}(\theta_{s+1})$.
	
	Firstly, if  we solve $p_{b}(\theta_s)\le p_{b}(\theta_{s+1})$, we get 
	\begin{align}b(\theta_{s}-\theta_{s+1})(\theta_{s+1}+\theta_{d}) &\le  (\theta_{s+1}-\theta_s )(\theta_{s+1}-\theta_d )(\theta_{s+1}+\theta_s+\theta_d) \nonumber\\
	b&\le-(\theta_s+\theta_{s+1}+\theta_d).  \label{b2}
	\end{align}

	%		 Thus the domain of $b$ is such that,  the smallest $b$ satisfies 
	%		$$\begin{cases}p_{b}(\theta_s)=p_{b}(\theta_{s-1}) &\text{ if } 1<s\le d-1\\
	%			p_{b}(\theta_s)<p_{b}(\theta_{s-1})& \text{ if } s=1,
	%		\end{cases}$$ and the largest $b$ satisfies $p_{b}(\theta_s)=p_{b}(\theta_{s+1})$ for all $1\le s\le d-1.$
	%		$$\begin{cases}p_{b}(\theta_s)=p_{b}(\theta_{s+1}) &\text{ if } 1\le s< d-1\\
	%		p_{b}(\theta_{d-1})=p_{b}(\theta_{d})& \text{ if } s=d-1,
	%		\end{cases}$$
	Secondly, if $1<s\le d-1$, then solving $p_{b}(\theta_s)\le p_{b}(\theta_{s-1})$ gives
	\begin{align}b(\theta_{s-1}-\theta_{s})(\theta_{d}-\theta_{s-1}) &\le  (\theta_{s-1}-\theta_s )(\theta_{s-1}-\theta_d )(\theta_{s-1}+\theta_s+\theta_d) \nonumber\\
	b&\ge-(\theta_s+\theta_{s-1}+\theta_d),  \label{b1}
	\end{align} and if $s=1$, then in order to satisfy the condition $\lambda(p_{b})=p_b(\theta_1)<p_{b}(\theta_{0})$,  $b$ must be such that $b>-(\theta_1+\theta_{0}+\theta_d)$. 
	Hence, if $\tau$ is an eigenvalue, $\theta_s$ of $G$, then
	\begin{align}\label{bbound}
	-(\theta_s+\theta_{s-1}+\theta_d) \le b\le -(\theta_s+\theta_{s+1}+\theta_d),
	\end{align}
	with a strict lower inequality when $s=1$.
	
	On the other hand, if $\tau=-\frac{\theta_{0}^2+\theta_{0}\theta_d-\triangle}{\theta_{0}(\theta_d+1)}$ is not an eigenvalue of $G$, then we can choose $\theta_s$ to be either the largest eigenvalue less than $\tau$ or the least eigenvalue greater than $\tau$.  
	\begin{itemize}
		\item If we choose $\theta_s$ to be the largest eigenvalue less than $\tau$, then $\varPhi(b)$ is an increasing function with a vertical asymptote at $b=-(\theta_{0}+\theta_s+\theta_d)$. We must pick the smallest $b$ that ensures $\lambda(p_{b})$  still occurs at $\theta_s$. Thus, $b$ must satisfy $p_{b}(\theta_s)=p_{b}(\theta_{s-1})$ %(since decreasing $b$ causes the minimum of $p_{b}(x)$ to transition from occurring at $\theta_s$ to $\theta_{s-1}$)
		yielding $b=-(\theta_1+\theta_{0}+\theta_d) $ if $1<s\le d-1$, and if $s=1$, then $b<-(\theta_1+\theta_{0}+\theta_d).$
		\item On the other hand, if we choose $\theta_s$ to be the least eigenvalue greater than $\tau$, then $\varPhi(b)$ is a decreasing function with a vertical asymptote at $b=-(\theta_{0}+\theta_s+\theta_d)$. We must pick the largest $b$ that ensures $\lambda(p)$ still occurs at $\theta_s$. Thus, $b$ must satisfy $p_{b}(\theta_s)=p_{b}(\theta_{s+1})$ yielding $b=-(\theta_s+\theta_{s+1}+\theta_d)$.
	\end{itemize}   
	
	Either of the above choices of $\theta_s$ will yield the same optimal bound since in each case $\lambda(p)=p(\theta_d)=p(\theta_s)$. Thus, we get a family of optimal polynomials. So, for simplicity, we  choose $\theta_s$ to be the least eigenvalue greater or equal to $\tau$ and then take $b=-(\theta_s+\theta_{s+1}+\theta_d)$. Now substituting $b$ into (\ref{cc}), we obtain $c=\theta_d\theta_s+\theta_d\theta_{s+1}+\theta_s\theta_{s+1}$. Finally, substituting $t=s$ or $t=d$ into (\ref{keyg}), and using $b$ and $c$ above, we obtain the desired result in (\ref{Gmaink-3}).
	
	Now the quotient matrix $B$ of $p(A)$ resulting from the proof of Theorem~\ref{t1} in \cite{abiad2019k} is
	\begin{align}
	B=\left[ 
	\begin{array}{cc}
	\frac{1}{r}\sum_{u\in U}(p(A))_{uu} &  p(\lambda_1)-\frac{1}{r}\sum_{u\in U}(p(A))_{uu}\\ \\
	\frac{rp(\lambda_1)-\sum_{u\in U}(p(A))_{uu} }{n-r}  &  p(\lambda_1)- \frac{rp(\lambda_1)-\sum_{u\in U}(p(A))_{uu} }{n-r}
	\end{array}
	\right], \label{Bmat}
	\end{align} with corresponding eigenvalues $\mu_1\ge \mu_2$,  and by interlacing \begin{align}\label{bIn}
	\lambda(p)\le \mu_2\le W(p)-\frac{rp(\lambda_1)-rW(p)}{n-r}.
	\end{align}
	So, if equality holds in (\ref{Gmaink-3}), then from (\ref{bIn}) we conclude that $\mu_2= \lambda(p)$ and, since $\mu_1 = p(\lambda_1)$, the interlacing is tight and
	the partition of $p(A)$ is regular (or equitable). Also, to derive its quotient matrix $B$ given in (\ref{B3mat}), we use (\ref{Bmat}) with the optimal polynomial $$p(x)=x^3-(\theta_s+\theta_{s+1}+\theta_d)x^2+(\theta_d\theta_s+\theta_d\theta_{s+1}+\theta_s\theta_{s+1})x$$ and the bound on $\alpha_3=r$ in (\ref{Gmaink-3}).
	Note that $W(p)=\frac{1}{r}\sum_{u\in U}(p(A))_{uu}$. Thus, $$B_{11}=W(p)=\triangle+b\theta_{0}=\triangle-(\theta_s+\theta_{s+1}+\theta_d)\delta.$$ We also have that $B_{12}=p(\delta)-W(p)$, so $$B_{12}=\delta^3-(\theta_s+\theta_{s+1}+\theta_d)\delta^2 + (\theta_d\theta_s+\theta_d\theta_{s+1}+\theta_s\theta_{s+1}+\theta_s+\theta_{s+1}+\theta_d)\delta-\triangle.$$ We have $$B_{21}=\frac{r}{n-r}(p(\delta)-W(p))=\triangle-(\theta_s+\theta_{s+1}+\theta_d)\delta-\theta_s\theta_{s+1}\theta_d.$$ Finally, $B_{22}=p(\delta)-B_{21}$ so
\begin{gather*}
B_{22}=\delta^3-(\theta_s+\theta_{s+1}+\theta_d)\delta^2 + (\theta_d\theta_s+\theta_d\theta_{s+1}+\theta_s\theta_{s+1}+\theta_s+\theta_{s+1}+\theta_d)\delta+\theta_s\theta_{s+1}\theta_d-\triangle.
	\end{gather*}
	
\end{proof}

\begin{remk}A polynomial $p$ satisfying $p(\theta_0)> \lambda(p)$  is optimal for $k=3$ if and only if  $$\lambda(p)=p(\theta_d)=p(\theta_s)=p(\theta_{s+1}),$$
	where
	$\theta_s$ is the least eigenvalue such that $\theta_s\ge -\frac{\theta_{0}^2+\theta_{0}\theta_d-\triangle}{\theta_{0}(\theta_d+1)}.$
\end{remk}

Observe that, if $G$ has girth greater than $3$, then $\triangle=0$, and the bound in Theorem~\ref{Gmaink3} becomes as follows.
\begin{corollary}\label{Gmaink3ii} Let $G$ be a $\delta$-regular graph with $n$ vertices, with girth greater than $3$, and distinct adjacency eigenvalues $\delta=\theta_0>\theta_1>\dots> \theta_d$, with $d\geq3$. Let $\theta_s$ be the least eigenvalue such that $\theta_s\ge -\frac{\theta_d+\theta_{0}}{\theta_d+1}$. %Let $b=-(\theta_s+\theta_{s+1}+\theta_d)$
	%and  $c=\theta_d\theta_s+\theta_d\theta_{s+1}+\theta_s\theta_{s+1}$. 
	Then, 
	\begin{align}\label{Gmaink-3i}
	\alpha_3%&\le n\frac{b\theta_0-\theta_s\theta_{s+1}\theta_d}{(\theta_{0}^3-\theta_d^3)+b(\theta_{0}^2-\theta_d^2)+c(\theta_{0}-\theta_{d}) }\\
	&\le n\frac{-\theta_{0}(\theta_s+\theta_{s+1}+\theta_d)-\theta_{s}\theta_{s+1}\theta_{d} }{(\theta_{0}-\theta_s)(\theta_{0}-\theta_{s+1})(\theta_{0}-\theta_{d})}.
	\end{align} 	
\end{corollary}

We use the result in Theorem~\ref{Gmaink3} to provide an upper bound for the diameter of a $\delta$-regular graph $G$.
\begin{corollary}
	Let $G$ be a $\delta$-regular graph with $n$ vertices, adjacency matrix $A$, and distinct eigenvalues $\delta=\theta_0>\theta_1>\dots>\theta_s>\theta_{s+1}>\dots, > \theta_d$, with $d\geq3$. Let $\theta_s$ be the least eigenvalue such that $\theta_s\ge -\frac{\theta_{0}^2+\theta_{0}\theta_d-\triangle}{\theta_{0}(\theta_d+1)}$, where $\triangle=\text{max}_{u\in V}\{(A^3)_{uu}\}$, twice the largest number of triangles on any vertex in $G$.
	%	Let $b=-(\theta_s+\theta_{s+1}+\theta_d)$
	%	and  $c=\theta_d\theta_s+\theta_d\theta_{s+1}+\theta_s\theta_{s+1}$. 
	If 
	\begin{align}\label{cor1}
	%n\frac{\triangle-\theta_{d}^3 +b(\theta_{0}-\theta_d^2)-c\theta_{d} }{(\theta_{0}^3-\theta_d^3)+b(\theta_{0}^2-\theta_d^2)+c(\theta_{0}-\theta_{d}) }<2,
	n\frac{\triangle-\theta_{0}(\theta_s+\theta_{s+1}+\theta_d)-\theta_{s}\theta_{s+1}\theta_{d} }{(\theta_{0}-\theta_s)(\theta_{0}-\theta_{s+1})(\theta_{0}-\theta_{d})}<2,
	\end{align} 
	%or for bipartite graphs, if \begin{align}\label{cor2}
	%	n\frac{(\theta_{0}^2+c) +b(1-\theta_0) }{2(\theta_{0}^2+c) }<2,
	%\end{align} 	
	
	then $G$ has diameter at most $3$.
\end{corollary}

\section{Some Applications}

To illustrate Theorem~\ref{Gmaink3}, we first consider an infinite family of graphs where the bound is tight.
\begin{defin} Let $G$ be a connected graph. 
	For a vertex $u\in V(G)$, we define $V_i(u)$ to be the set of vertices at distance $i$ from
	$u$. Then, $G$ is a \textit{distance-regular graph}  if  $|V_i(u) \cap  V_j(v)|$ depends only on 	the distance between vertices $u$ and $v$.
\end{defin} Thus, given any two vertices $u$ and $v$ at distance $k$ in a  distance-regular graph $G$, the number of vertices at distance $i$ from $u$ and distance $j$ from $v$ is determined by $k, i, j$. An {\em intersection array} of graph $G$ of diameter $d$ is a sequence of integers $\{b_0,b_1,\dots, d_{d-1},c_1,\dots,c_d\}$ such that for all $1\le k\le d$, $b_k$ is the number of neighbours of $u$ at distance $k+1$ from $v$ and $c_k$ is the number of neighbours of $u$ at distance $k-1$ from $v$. For a detailed treatment of distance-regular graphs, see Brouwer, Cohen and Neumaier \cite{brouwer1989}. 

A graph $G$ of diameter $d$ is {\em antipodal} if there exists a partition of the vertex set into classes with the property that any two distinct
vertices in the same class are at distance $d$, while two vertices in different
classes are at distance less than $d$.
\subsection{Antipodal bipartite distance-regular graphs} %\label{antipodal}
	Let $G$ be an antipodal bipartite distance-regular graph, with degree  $d$ and diameter $3$. These graphs have $n=2(d+1)$ vertices (Brouwer, Cohen and Neumaier \cite{brouwer1989}), intersection array $\{d,d-1,1;1,d-1,d\}$, and distinct eigenvalues \begin{align}
	\theta_{0}=d, \, \theta_{1}=1, \, \theta_{2}=-1, \, \theta_{3}=-d.
	\end{align} 
	As $G$ is bipartite, $\theta_{s}=1$. Thus, by Theorem \ref{Gmaink3}, we have 
	\begin{align*}
	\alpha_3 &\le 2(d+1)\frac{d-(1-1)+(1)(-1)}{2(d-1^)(d+1)}\\
	&= 2(d+1)\frac{(d-1)}{2(d-1)(d+1)}=1\\
	\end{align*}
	as expected since $G$ has diameter $3$. Also, since $b=d$ and $c=-1$, the polynomial that gives this bound is $p(x)=x^3+dx^2-x$.
	
	If the adjacency matrix of the graph is $\mathbf{A}$, then we have that the matrix $p(\mathbf{A})=\mathbf{A}^3+d\mathbf{A}^2-\mathbf{A}$ has a regular partition with the quotient matrix
	\begin{align*}
	B=\left[\begin{array}{cc}
	d^2 & 2d^3-d^2-d\\d^2+d& 2d^3-d^2-2d
	\end{array}\right].
	\end{align*}

\subsection{Hamming Graphs}
Let $d\ge 3$ and $q\ge 2$ be integers, and let $Q$ be a set of $q$ elements. The Hamming graph $H(d,q)$ is a $d(q-1)$-regular graph with vertex set $Q^d$, consisting of sequences of length $d$ from $Q$ and two sequences (vertices) are adjacent if they differ in just one position. A 3-independent set in $H(d,q)$ thus, consists of sequences that differ in at least 4 positions. $H(d,q)$ has $q^d$ vertices and eigenvalues $d(q-1)-qi$ for $i=0,1,\dots, d$ with respective multiplicities $\binom{d}{i}(q-1)^i$. See Brouwer and Haemers \cite{BrouwerHaemer}. Also,  $\triangle =2d\binom{q-1}{2}=d(q-1)(q-2)$, $\theta_0=d(q-1)$,  $\theta_{d}=-d$ and $s$ is the largest index such that \thesisonly{ \begin{align*}
    \theta_s &\ge -\frac{\theta_{0}^2+\theta_{0}\theta_d-\triangle}{\theta_{0}(\theta_d+1)}\\
    &= -\frac{d^2(q-1)^2-d(q-1)d-d(q-1)(q-2)}{d(q-1)(-d+1)}\\
    &= -\frac{d(q-1)-d-(q-2)}{(-d+1)}=\frac{(d-1)(q-2)
    }{d-1} =q-2.
\end{align*} Thus, $s$ is the largest index such that} $\theta_{s}\ge q-2$. Now let $b=d\mod{q}$, then solving $\max \{s: d(q-1)-sq\ge q-2\}$,  \thesisonly{that is, \begin{align*}
    d(q-1)-sq &\ge q-2\\
    d-1-\frac{d-2}{q}&\ge s\\
    d-1-\Big\lceil{\frac{d-2}{q}}\Big\rceil &= s.
\end{align*}
Now, $d=aq+b$ for some $a\in \mathbb{Z}$, so 
\begin{align*}
    s &=d-1-\Big\lceil{\frac{d-2}{q}}\Big\rceil\\
    &= d-1-\Big\lceil{\frac{aq+b-2}{q}}\Big\rceil\\
    &= d-1-a - \Big\lceil{\frac{b-2}{q}}\Big\rceil.\\
\end{align*}
Now if $q=2$, we have 
\begin{align*}
    s = d-1-a - \Big\lceil{\frac{b-2}{q}}\Big\rceil=d-a-\Big\lceil{\frac{b}{2}}\Big\rceil= \begin{cases}
        d-a \text{ if } b=0\\
        d-a-1 \text{ if } b=1.
    \end{cases}
\end{align*}
Similarly, if $q>2$, we have 
\begin{align*}
    s = d-1-a - \Big\lceil{\frac{b-2}{q}}\Big\rceil= \begin{cases}
        d-a-1 \text{ if } 0\le b \le 2\\
        d-a-2 \text{ if } 2< b < q.
    \end{cases}
\end{align*}
Now, substituting the index $s$,} we have  $\theta_{s}=0$ and $\theta_{s+1}= -2$ if $q=2$ and $b=0$, otherwise 
\[ \theta_s= 
\begin{cases}
    -b+q \text{ if } 0\le b\le 2\\
    -b+2q \text{ if } 2< b <q,
\end{cases}\]
and 
\[ \theta_{s+1}= 
\begin{cases}
    -b \text{ if } 0\le b\le 2\\
    -b+q \text{ if } 2< b <q.
\end{cases}\]
We deduce the following from Theorem~\ref{Gmaink3}.
\begin{proposition}\label{genHam}
    Let $b=d\mod{q}$. The $3$-independence number of the Hamming graph $H(d,q)$ satisfies the following.
    %\begin{align}
% \alpha_3(H(d,q))  \le
% \begin{cases}\label{Genhamm}
%     q^{d-1} \frac{q(d+b-2)+(b-1)^2+(1-d)}{(dq-d+b)(dq-d+b-q)} \text{ if } 0\le b \le 2\\
%     \\
%     q^{d-1} \frac{q(d-b)+(b-1)^2+(1-d)}{(dq-d+b-2q)(dq-d+b-q)} \text{ if } 2<b <q.
% \end{cases}\end{align}
%\begin{numcases}{\alpha_3(H(d,q))  \leq }  q^{d-1} \frac{q(d+b-2)+(b-1)^2+(1-d)}{(dq-d+b)(dq-d+b-q)} \text{ if } 0\le b \le 2 \label{Genhamm1}
%  \\ 
 %  \nonumber\\
  %  q^{d-1} \frac{q(d-b)+(b-1)^2+(1-d)}{(dq-d+b-2q)(dq-d+b-q)} \text{ if } 2<b <q. \label{Genhamm2}
%\end{numcases}
 
\begin{numcases}{\alpha_3(H(d,q))  \leq }q^{d-1} \frac{d-2}{d(d(q-1)-q)} \text{ if } b=0 \label{GenBB0}\\
	\nonumber\\
	q^{d-1}\frac{1}{d(q-1)+1} \text{ if } b=1 \label{GenBB}\\
	\nonumber\\
	q^{d-1}\frac{1}{d(q-1)-q+2} \text{ if } b=2 \label{Genhamm1}
		\\ 
		\nonumber\\
		q^{d-1} \frac{q(d-b)+(b-1)^2+(1-d)}{(dq-d+b-2q)(dq-d+b-q)} \text{ if } 2<b <q. \label{Genhamm2}
	\end{numcases}

\end{proposition}
 %To investigate tightness, let us first consider when we have integer bounds.  We can simplify (\ref{Genhamm1}) further when we consider the cases $b=0,1,2$. %Note that $b=0=2$ for $q=2$.
%\begin{corollary}\label{genb}
 %   Let $b=d\mod{q}$. The $3$-independence number of the Hamming graph $H(d,q)$ satisfies the following.
%    \begin{align}
% \alpha_3(H(d,q))  \le \begin{cases}\label{GenBB}
%     q^{d-1} \frac{d-2}{d(d(q-1)-q)} \text{ if } b=0\\
%    \\
 %     \frac{q^{d-1}}{d(q-1)+1} \text{ if } b=1\\
 %     \\
   %   \frac{q^{d-1}}{d(q-1)-q+2} \text{ if } b=2.
% \end{cases}  \end{align}
% \end{corollary}

It is easy to see that for $d=q$, the bound on $H(d,q)$ is an integer. Indeed, $\alpha_3(H(d,d))= \alpha_3(H(q,q))\le q^{q-3}\in \mathbb{Z}^+$. More generally, we have an integer bound on $H(d,q)$  for $3\le d\le q+2$ with $q\ge 2$. % and $\alpha_3(H(q,q))\le q^{d-3}\in \mathbb{Z}^+$.
Observe that $b=d$ for $d\in [3,q-1]$, so substituting $b=d$ into Equation (\ref{Genhamm2}) of Proposition~\ref{genHam}, we have an integer bound $q^{d-3}$. Moreover, we have $b=0$, $b=1$ and $b=2$ for $H(q,q)$, $H(q+1,q)$ and $H(q+2,q)$ respectively, hence, Equations (\ref{GenBB0}-\ref{Genhamm1}) each yields an integer bound $q^{d-3}$ for each case.  We summarize this result below.
\begin{corollary}\label{intbound}
    Given  $3\le d\le q+2$,  the right side of Equations (\ref{GenBB0}-\ref{Genhamm2}) simplifies to $q^{d-3}$, therefore
    \begin{align}
\alpha_3(H(d,q))  \le
    q^{d-3}.
 \end{align}
 \end{corollary}
 The Singleton bound~\cite{stinson2007combinatorial} (on the number of codewords in a code of length $d$ and minimum distance 4  ) gives $\alpha_3\le q^{d-3}$ for all $d$ and $q$, and when $d> q+2$ the ratios of the right sides of the bounds in Proposition~\ref{genHam} to the Singleton bound approach 0 as $d$ approaches infinity.
%\thesisonly{Observe that the bound in Corollary~\ref{intbound} coincides with the Singleton Bound~\cite{} on the number of codewords in a code of minimum distance 4. Tight bounds will therefore provide Maximum Distance Separable (MDS) codes.}

Now, when $b=1$, bound (\ref{GenBB}) is  an integer if  $d=\frac{q^r-1}{q-1}$ with $r\in \mathbb{Z}^+ $. To see this, observe that the denominator of the bound is \begin{align*}
    d(q-1)+1=\frac{q^r-1}{q-1}(q-1)+1=q^r.
\end{align*}
Thus, \begin{align*}
    \alpha_3(H(d,q))\le \frac{q^{d-1}}{d(q-1)+1}=\frac{q^{\frac{q^r-1}{q-1}-1}}{q^r}=q^{\frac{q^r-q}{q-1}-r}\in \mathbb{Z}^+.
\end{align*}
Similarly, when $b=2$, bound (\ref{Genhamm1}) is an integer if  $d=\frac{q^r+(q-2)}{q-1}$ with $r\in \mathbb{Z}^+$. Observe that the denominator of the bound is \begin{align*}
    d(q-1)-q+2= \frac{q^r+(q-2)}{q-1}(q-1)-q+2 =q^r.
\end{align*}
Thus, \begin{align*}
    \alpha_3(H(d,q))\le \frac{q^{d-1}}{d(q-1)-q+2}=\frac{q^{\frac{q^r-1}{q-1}}}{q^r}=q^{\frac{q^r-1}{q-1}-r}\in \mathbb{Z}^+.
\end{align*}
Hence, one can study tightness at these values of $d$ and $q$ that lead to integer bounds.
Let us now consider some special Hamming graphs and investigate their tightness.
\subsubsection{The Hamming graph $H(3,q)$}
Consider the Hamming graph $H(3,q)$. It has $q^3$ vertices and the distinct eigenvalues $3q-3,2q-3$, $q-3$ and $-3$. Bang, van Dam and Koolen \cite{MR2455522} showed $H (3,q)$  is uniquely determined by its spectrum for $q\ge 36$. % We have $\triangle =3(q-1)(q-2)$, $\theta_s=2q-3$, and $\theta_{s+1}=q-3$. 
Now, $3\mod{q}$ is $0,1$ or $3$ for $q\geq 2$. Thus, Equations (\ref{GenBB0}), (\ref{GenBB}) and (\ref{Genhamm2}) give exactly $1$  for $\alpha_3(H(3,q))$,  %$3$-independence number of the Hamming graph $H(3,q)$ with $q> 3$, is  
%\begin{align*}     \alpha_3(H(3,q))  &\le q^3 \frac{3(q-1)(q-2)-(3q-3)(3q-9)+3(2q-3)(q-3)}{(q)(2q)(3q)}\\     &= 1 \end{align*}
as expected, since $H(3,q)$ has diameter 3.

\subsubsection{The Hamming graph $H(4,q)$}
Consider the Hamming graph $H(4,q)$. Using Corollary~\ref{intbound}, we have that $\alpha_3(H(4,q))  \le q$.  This is a tight bound. To see this, consider a $3$-independent set of $H(4,q)$ given as \{1111, 2222, 3333, \dots, qqqq\}, which has cardinality $q$.

\subsubsection{The Hamming graph $H(5,q)$}
Corollary~\ref{intbound} shows that $\alpha_3(H(5,q))  \le q^2$.\thesisonly{Recall a \emph{Latin square $L$} of order $q$ is an $q \times q$ array in which $q$ distinct symbols $\{1,2,3,\dots, q\}$ are arranged so that
each symbol occurs exactly once in each row and column. Two Latin squares $L_1$ and $L_2$ of order $q$ are said to be \emph{orthogonal} if when superimposed, each of the possible $q^2$ pairs occurs exactly once. A set of Latin squares is said to be \emph{mutually orthogonal} if  the Latin squares are pairwise orthogonal. This set is termed Mutually Orthogonal Latin Squares (MOLS). Suppose $\ell$ is the number of MOLS of size $q$ that exists, then since each MOLS is a $q\times q$ matrix, we can create $q^2$ codes of length $\ell +2$ by creating an orthogonal array (OA) from the MOLS (see Stinson~\cite{stinson2007combinatorial}).} If there are $\ell$  Mutually Orthogonal Latin Squares (MOLS) of size $q$, then since each Latin square is a $q\times q$ matrix, we have an orthogonal array of $q^2$ rows  and $\ell +2$ columns (see Stinson~\cite{stinson2007combinatorial}). The rows of this array form a code with minimum distance $\ell+1$. Thus, MOLS codes give constructions for a 3-independent set of  $H(5,q)$ that meet our bound when there are at least $\ell=3$ MOLS of size $q$. For instance, for $q$ a prime power, $\ell=q-1$. From Hanani~\cite{hanani1970number}, Wilson~\cite{wilson1974concerning} and  the survey work of Colbourn and  Dinitz~\cite{colbourn2001mutually}, we know the lower bounds on $\ell$ for $q\ge 3$. Thus, since $\ell\ge 3$ for $q\ge 3$ except  for $q=3,6$ and possibly for $q=10$, we have a 3-independent set for $H(5,q)$ which meets our bound.
 
\subsubsection{The Hamming graph $H(d,3)$}
The Hamming graph $H(d,3)$ has distinct eigenvalues $2d-3i$ for $i=0,1,2,\dots, d.$ It has $\triangle=2d$,  $\theta_s=3-b$ and $\theta_{s+1}=-b$, where $b=d\mod{3}$. We have the following bound for $H(d,3)$. 
\begin{corollary}
    Let $b=d\mod{3}$. The $3$-independence number of the Hamming graph $H(d,3)$ satisfies the following.
    \begin{align}
\alpha_3(H(d,3))  \le
    3^{d-1} \frac{2d+b^2+b-4}{(2d+b)(2d+b-3)}.
 \end{align}
Alternatively,
\begin{align}
    \alpha_3(H(d,3))  \le
\begin{cases}\label{hamd3}
    3^{d-1}\frac{d-2}{d(2d-3)} \text{ if } b=0\\ \\
    \frac{3^{d-1}}{2d+1} \text{ if } b=1\\ \\
    \frac{3^{d-1}}{2d-1} \text{ if } b=2.
\end{cases}
\end{align}

\end{corollary}
It is easily verified that the bound in (\ref{hamd3}) is tight for $d=3$. Note that the bound is an integer when $d=\frac{3^r-1}{2}$ for $b=1$, and when $d=\frac{3^r+1}{2}$ for $b=2$ for all $r\in \mathbb{Z}^{+}$. %Thus, any possibility of tightness must occur at these values of $d$.  %We verify(ied) using GAP tightness for $d=1,4,13,40,121$ and $d=2,5,14,41,122$.
It is also an integer when $d=6$. Using GAP, we confirm $H(6,3)=18$, meaning it is tight. It's however not tight for $H(5,3)$, as bound is 9 but actual value is 6. Thus, the bound (\ref{hamd3}) is not tight in general.
\subsubsection{ The Hamming graph $H(d,2)$}
Consider the Hamming graph %{\sc Hamming graph} 
$H(d, 2)$ also known as the hypercube $Q_d$ or the $d$-cube. The vertices of this graph are binary
sequences (hence forth, sequences) of length $d$ and two sequences (vertices) are adjacent if they differ in just one position. Thus, it has $2^d$ vertices. Suppose the matrix $\mathbf{A}$ is its adjacency matrix. Then, its eigenvalues are the integers
$d-2i$, for $i = 0,...,d$
with respective multiplicities $\binom{d}{i}$. See Brouwer and Haemers~\cite{BrouwerHaemer}  for further details. It is bipartite and so $\triangle=0$ and $\theta_{0}=-\theta_d=d$. If its diameter $d$ is odd, then $\theta_s=1$ and $\theta_{s+1}=-1$. Moreover, if $d$ is even, then $\theta_s=0$ and $\theta_{s+1}=-2$. Thus, by Theorem~\ref{Gmaink3}, we have the following. %that, for $d$ odd
% \begin{align}
% \alpha_3(H(d, 2)) \le  2^d \frac{d-(1-1)-1}{2(d-1)(d+1)}=\frac{2^{d-1}}{d+1} \label{hamm1}
% \end{align}
% and for $d$ even
% \begin{align}
% \alpha_3(H(d, 2)) \le  2^d \frac{d-(0-2)}{2(d)(d+2)}=\frac{2^{d-1}}{d}.\label{hamm2}
% \end{align}

\begin{corollary} The $3$-independence number of the Hamming graph $H(d,2)$ satisfies the following.
%\begin{equation}
    \begin{numcases}{\alpha_3(H(d, 2)) \leq }
 2^d \frac{d-(1-1)-1}{2(d-1)(d+1)}=\frac{2^{d-1}}{d+1} \text{ if } d \text{ is odd } \label{hamm1}
  \\ 
   \nonumber\\
    2^d \frac{d-(0-2)}{2(d)(d+2)}=\frac{2^{d-1}}{d} \text{ if } d \text{ is even.} \label{hamm2}
\end{numcases}
%\end{equation}
 
\end{corollary}

%\end{align}
These bounds coincide with the bounds of Fiol in \cite{fiol2}.
  
%There is some evidence of tightness for the case $d=2^i$ for all $i\in \mathbb{Z}^{+}$, but we have not pursued a proof or otherwise.
%\textcolor{red}
{In what follows, we construct a $3$-independent set and thus, a lower bound for $\alpha_3(H(d,2))$. We then show that bound (\ref{hamm2}) is tight for some family of $d$, specifically $d=2^r$, while bound (\ref{hamm1}) is tight for $d=2^r-1$,  for all $r\in \mathbb{Z}^{+}$.}

\subsubsection{A $3$-independent set for $H(d,2)$ }

 We first present a general result.
\begin{proposition}\label{constG}
    Let $d_1\le d_2$. Suppose $\alpha_3(H(d_1,2))=\ell_1$ and $\alpha_3(H(d_2,2))=\ell_2$, then $\alpha_3(H(d_1+d_2,2))\ge \ell_1\ell_2d_1.$
\end{proposition}
\begin{proof}
Given $d_1\le d_2$, take a $3$-independent set $U$ of size $\ell_1$ in $H(d_1,2)$  and a $3$-independent set $V$ of size $\ell_2$ in $H(d_2,2)$.
We construct a $3$-independent set of size $\ell_1\ell_2d_1$ in $H(d_1+d_2,2)$.
 Let $W=U\times V$ be the cartesian product of $U$ and $V$. Thus, the size of $W$ is $\ell_1\ell_2$. 
% $\frac{2^{d-1}}{d}$ for $d=2^r$. 
%Take the Cartesian product  of $W$ with itself to get a set $W^2$ in $H(2d,2)$ of size $  \left(\alpha_3(H(d,2))\right)^2$.
Any pair of sequences $x_1y_1$ and $x_2y_2$  in $W$ with $x_1,x_2\in U$ and $y_1,y_2\in V$, are such that $x_1\ne x_2$ or $y_1\ne y_2$ or both. So, $x_1y_1$ and $x_2y_2$ are $3$-independent in $W$ since $x_1$ and $ x_2$ are $3$-independent in $U$  or $y_1$ and $ y_2$ are $3$-independent in $V$. Hence, $W$ is a $3$-independent set.  % Now, for each sequence in $W^2$, and for $j=1,2,..,d$ flip the columns $j$ and $d+j$.
For each $j=1,2,..,d_1$, let $W_j$ be the  set of sequences resulting from  adding $1\mod 2$  to the  $j$ and $(d_1+j)$th elements of each sequence in $W$. (Hence forth, we shall say the $i$th column of a sequence is flipped when $1\mod 2$ is added to the $i$th element of the sequence).
 This results in a set $\displaystyle A=\bigcup_{j=1}^{d_1} {W_j}$ of size $\ell_1\ell_2 d_1$. We claim set $A$ is $3$-independent. For simplicity, assume the sequences in a set are rows of a matrix representing the set.
 First of all, any pair of sequences from the same $W_j\in A$ is a result of flipping the same columns of sequences in $W$, so their Hamming distance remains unchanged. Also, any pair of sequences each from $W_i$ and $W_j$, with $i\ne j$, either resulted from the same sequence or two distinct sequences of $W$. If the pair are from the same sequence, then we have $4$ distinct columns of the sequence flipped, two for each pair. Hence, the pair of sequences have a Hamming distance exactly $4$. On the other hand, if the pair is from two distinct sequences, then at least four out of their Hamming distance came from either the first $d_1$ columns or the last $d_2$ columns, so the flip reduced the Hamming distance by at most two in that half but also gained two in the other half. Hence, the Hamming distance remains at least $4$.
 Thus, we have obtained a $3$-independent set $A$ for $H(d_1+d_2,2)$ of size $ \ell_1\ell_2d_1$ proving the claim.
\end{proof}
An immediate consequence of Proposition \ref{constG} is as follows.
 If $d=d_1= d_2$ and $\alpha_3(H(d,2)) =\ell$, then we have obtained a $3$-independent set $A$ for $H(2d,2)$ of size $d\ell^2$.
 We now show by induction that it  coincides with the upper bound (\ref{hamm2}) when $d=2^r$, for all $r\in \mathbb{Z}^{+}$, and thus, tight. It is clear that for $r=1$, that is, $d=2$, we have just one element in the $3$-independent set, say $00$. Suppose we have a $3$-independent set 
 in $H(d,2)$ of size  $\displaystyle
\ell=\frac{2^{d-1}}{d} $ for $d=2^r$,  then the construction gives  $\displaystyle d\ell^2= \frac{2^{2d-1}}{2d}$ for $2^{r+1}=2d$. But $\displaystyle\frac{2^{2d-1}}{2d}=\alpha_3(H(2d, 2))$, using bound
(\ref{hamm2}).  % But \[\frac{2^{2d-1}}{2d}=\alpha_3(H(2d, 2)).\]
We have thus, shown that the bound (\ref{hamm2}) is tight when $d=2^r$, and also proved the following result.
\begin{corollary}
    If $\alpha_3(H(d,2))=\ell$, then $\displaystyle\alpha_3(H(2d,2))\ge d\ell^2$ with equality if $d=2^r$, for all $r\in \mathbb{Z}^{+}$.
\end{corollary}

In the case where equality is attained, Theorem~\ref{Gmaink3} tells us more. We have the following.
\begin{corollary}
    Let $\mathbf{A}$ be the adjacency matrix of the Hamming graph $H(d,2)$ with $d=2^r$, for all $r\in \mathbb{Z}^{+}$. Then, the matrix $\mathbf{A}^3+(2+d)\mathbf{A}^2+2d\mathbf{A}$ has a regular partition
    (with a set of $\frac{2^{d-1}}{d}$ $3$-independent vertices and their complement) with the quotient matrix
	\begin{align*}
	\mathbf{B}=\left[\begin{array}{cc}
	d^2+2d & d^3+2d^2\\
        d^2+2d& 2d^3+3d^2-2d
	\end{array}\right].
	\end{align*}
\end{corollary}
\thesisonly{Recall the elements of the quotient matrix $\mathbf{B}$ are
\begin{align*}
   B_{11}&=0-(-2-d)d=d^2+2d\\
   B_{12}&=d^3-(-2-d)d+(2d-2-d)-0=d^3+2d\\
   B_{21}&=0-(-2-d)d-0=d^2+2d\\
   B_{22}&=d^3-(-2-d)d+(2d-2-d)d+0-0=3d^2-2d.
\end{align*}
}
% As corollary, we can construct $3$-independent sets and thus lower bounds of $\alpha_3(H(l,2))$ for $d< l<2d$. We start with a $3$-independent set in $H(2d,2)$ of size $\frac{2^{2d-1}}{2d}$ with $d=2^r$.  Without loss of generality, take a column that has as many $1$s as $0$s, and pick the sequences that contain the $1$s.  Delete these $1$s and we have found a set of sequences of length $2d-1$. It is easy to see this set is still $3$-independent. Moreover, there are half the size of the $3$-independent set in $H(2d,2)$ contained in this set. Hence $$\alpha_3(H(2d-1,2))\ge \frac{2^{2d-1}}{4d}.$$ Thus, recursively, we can construct a $3$-independent set for $H(l,2)$ for all $d< l<2d$. %  We can repeat this process to get $\alpha_3(H(2d-2,2))=\frac{2^{2d-1}}{8d}$, and also until $\alpha_3(H(d+1,2))$.  Thus, the $3$-independence number decreases by a factor of $2$ from $\alpha_3(H(2d,2))=\frac{2^{2d-1}}{2d}$ till $\alpha_3(H(d+1,2))$%=\frac{2^{d}}{d+1}$ \, with $d=2^r$, where $r\in \mathbb{Z}^{+}$.

Let $X$ be a set of cardinality $n$, consisting of sequences of length $d$. Create a matrix $\mathbf{X}_{n,b}$ from $X$ by making each of the $n$ sequences a row. We say $X$ is \emph{balanced} if each column of $\mathbf{X}_{n,b}$ has as many $1$s as $0$s. We show briefly that the constructed 3-independent set $A$ above is balanced for all $r>1$. Suppose we start with a balanced $3$-independent set $W$. Note that for $r=2$, the construction gives a $3$-independent set as $W=\{1010,0101\}$ which is balanced. Now, $W^2$ is also balanced since it is a concatenation of sequences of a balanced set. Moreover, $W^2_j$ is balanced since flipping the $j$ and $j+d$ columns of sequences of $W^2$ just swaps $0$s and $1$s in those columns. Thus, $A$ is balanced, being the union of balanced sets. 

Now, as corollary of the construction above, we can construct $3$-independent sets and thus, lower bounds of $\alpha_3(H(l,2))$ for $\frac{1}{2}d< l<d$. We start with a $3$-independent set  in $H(d,2)$ of size $\ell$ with $d=2^r$. Recall that this set is balanced for all $r\ge 2$.
Pick any column and take the sequences that contain a $1$ in that column.
 Delete these $1$s and we have found a set of sequences of length $d-1$. It is easy to see that this set is still $3$-independent, and thus, a $3$-independent set of size $\frac{1}{2}\ell$ in $H(d-1,2)$. We summarize the above in the following.
 \begin{proposition}\label{hamcom}
 The $3$-independence number of the Hamming graph  satisfies  \begin{align*}
\alpha_3(H(d-1, 2)) \ge \frac{1}{2}\alpha_3(H(d, 2)).
\end{align*} 
 \end{proposition} 
 
 Thus, recursively, we can construct a $3$-independent set for $H(l,2)$ for all $\frac{1}{2}d< l<d$.
 
 \begin{proposition}\label{hamconst} The $3$-independence number of the Hamming graph  satisfies %the following:
   \begin{align}
\alpha_3(H(d-i, 2)) \ge \frac{2^{d-(i+1)}}{d} \label{hamconst1}
\end{align}
with $0\le i\le \frac{1}{2}(d-2)$ and  $d=2^r$, for all $r\in \mathbb{Z}^{+}$.
\end{proposition}

 An immediate consequence of these is that bound (\ref{hamm1}) is tight for $d=2^r-1$ for all $r\in \mathbb{Z}^{+}$. To see this, we know that for $d=2^r$ we have  \begin{align*}
\alpha_3(H(d, 2)) = \frac{2^{d-1}}{d}.
\end{align*} So, from Proposition~\ref{hamcom}, we have \begin{align*}
\alpha_3(H(d-1, 2)) \ge \frac{1}{2} \alpha_3(H(d, 2))=\frac{1}{2}\frac{2^{d-1}}{d}=\frac{2^{d-2}}{d},
\end{align*} while from bound (\ref{hamm1}), we have \begin{align*}
\alpha_3(H(d-1, 2)) \le \frac{2^{(d-1)-1}}{(d-1)+1}=\frac{2^{d-2}}{d}.
\end{align*} Thus, \begin{align*}
\alpha_3(H(d-1, 2)) = \frac{2^{d-2}}{d},
\end{align*}  proving that bound (\ref{hamm1}) is tight for $d-1=2^r-1$ for all $r\in \mathbb{Z}^{+}$. That is, 
\begin{align*}
\alpha_3(H(d, 2)) = \frac{2^{d-1}}{d}
\end{align*}
for $d=2^r-1$ for all $r\in \mathbb{Z}^{+}$.

The following is immediate.
\begin{corollary}
    Let $\mathbf{A}$ be the adjacency matrix of the Hamming graph $H(d,2)$ with $d=2^r-1$, for all $r\in \mathbb{Z}^{+}$. Then, the matrix $\mathbf{A}^3+d\mathbf{A}^2-\mathbf{A}$ has a regular partition
    (with a set of $\frac{2^{d-1}}{d+1}$ $3$-independent vertices and their complement) with the quotient matrix
	\begin{align*}
	\mathbf{B}=\left[\begin{array}{cc}
	d^2 & 2d^3-d^2-d\\
        d^2-d & 2d^3-d^2-2d
	\end{array}\right].
	\end{align*}
\end{corollary}

\thesisonly{Recall the elements of the quotient matrix $\mathbf{A}$ are
\begin{align*}
   B_{11}&=0-(1-1-d)d=d^2\\
   B_{12}&=d^3-(-d)d+(-d+d-1-d)d-0=d^3+2d\\
   B_{21}&=0-(-2-d)d-0=d^2+2d\\
   B_{22}&=d^3-(-2-d)d+(2d-2-d)d+0-0=3d^2-2d.
\end{align*}
}

%See table~\ref{ncube}.
\begin{table}[h!]
	\caption{Bound for some Hamming graphs $H(d, 2)$} 
	\centering
	{\begin{tabular}{|c|c|c|c|}
			\hline 
			 $d$-cube & Theorem~\ref{Gmaink3} & Proposition~\ref{hamconst} & $\alpha_3$\\ 
			\hline
               2 &   1 & 1 & 1 \\ 
			\hline 
			3 &   1 & 1 & 1 \\ 
			\hline 
			4 & 2 &2  & 2\\ 
			\hline 
			5 & 2.6 &2 &2 \\ 
			\hline 
			6 & 5.3 &4 & 4  \\ 
			\hline 
			7 & 8 & 8& 8  \\ 
			\hline 
			8 & 16 & 16 & 16\\ 
   \hline 
			9 & 25.6 & 16 & 20\\ 
		%	\hline 
		%	J(16,8) & 12 & 4 & 44 &-\\ 
		%	\hline 
		%	J(18,9) & 11 & 3 & 117.54 &-\\
		%	\hline
		%	J(20,10) & 10 & 2 & 342.78 &-\\
			\hline
	\end{tabular} }\label{ncube}
\end{table}
 Observe from Table~\ref{ncube} that  bound~(\ref{hamconst1}) in Proposition~\ref{hamconst} is tight for all $d\le 8$. Observe also that $\alpha_3(H(8, 2))=16$ and $\alpha_3(H(9, 2))=20$, so by Proposition~\ref{constG} $\alpha_3(H(17, 2))\ge 2560$, while Proposition \ref{hamconst} gives $\alpha_3(H(17, 2))\ge\frac{2^{16}}{32}=2048$. Hence, Proposition~\ref{constG} can be stronger than Proposition~\ref{hamconst}.
We can also compare  the upper bounds (\ref{hamm1} and \ref{hamm2})  with the size of the $3$-independent sets constructed. 
Observe that the upper bounds (\ref{hamm1} and \ref{hamm2}) can be summarized as in the following proposition.
\begin{proposition}\label{ham12} The $3$-independence number of the Hamming graph  satisfies %the following:
   \begin{align}
\alpha_3(H(d-i, 2)) \le \frac{2^{d-(i+1)}}{d-i} \label{ham121}
\end{align}
with $0\le i\le \frac{1}{2}(d-2)$ and  $d=2^r$, for all $r\in \mathbb{Z}^{+}$.
\end{proposition}

%Now, the ratio of the bound (taking the right  side of Equation~(\ref{ham121})) to the size of the constructed $3$-independent set (taking the right side of Equation~(\ref{hamconst1})) is \begin{align*}
%\frac{d}{d-i}
%\end{align*}for both $H(d-i-1, 2)$ and $H(d-i, 2)$ with $d=2^r$, and for all even $i$ with $2\le i\le \frac{1}{2}(d-2)$. %The ratio is therefore asymptotically $2$ ($1$?)\footnote{Within each interval of $\frac{1}{2}d$ and $d$, the ratio varies from $1$ when $i=2$ to at most $2$ when $i = \frac{1}{2}(d-2)$. On the other hand, the ratio approaches $1$ when $d\longrightarrow \infty$. I'm not sure which one makes the most sense, or what we want to state there.}. 
%For $r$ fixed and for all even $i$ with $2\le i\le \frac{1}{2}(d-2)$, the ratio decreases as $d$ increases. Moreover, for $i = \frac{1}{2}(d-2)$ and $r$ increasing, the ratio asymptotically approaches $2$.

Now, the ratio of the bound (taking the right side of Equation~(\ref{ham121})) to the size of the constructed $3$-independent set (taking the right side of Equation~(\ref{hamconst1})) is \begin{align*}
\frac{d}{d-2j}
\end{align*}
for both $H(d-2j-1, 2)$ and $H(d-2j, 2)$ with $d=2^r$, and  $0\le j\le \frac{1}{4}(d-2)$. 
For $r$ fixed and $0\le j\le \frac{1}{4}(d-2)$, the ratio decreases as $d$ increases. Moreover, for $j= \frac{1}{4}(d-2)$ and $r$ increasing, the ratio asymptotically approaches $2$.

%Similarly, the difference between the bound (taking the right  side of Equation~(\ref{ham121})) and the size of the constructed $3$-independent set (taking the right side of Equation~(\ref{hamconst1})) is \begin{align*} \frac{i2^{d-(i+1)}}{d(d-i)} \end{align*} for $H(d-i, 2)$, while the difference between the bound and the size of the constructed set is \begin{align*} \frac{i2^{d-(i+2)}}{d(d-i)} \end{align*} for $H(d-i-1, 2)$ with $d=2^r$, and for all even $i$ with $2\le i\le \frac{1}{2}(d-2)$.

\subsection{Odd graphs}
Consider also the  Odd graph $O_\ell$. The vertices of $O_\ell$ correspond to the $(\ell-1)$-element subsets of a $(2\ell-1)$-element set, and two vertices are adjacent if and only if the corresponding subsets are disjoint. The Odd graph $O_\ell$ is an $\ell$-regular graph of order $n=\binom{2\ell-1}{\ell-1}=\frac{1}{2}\binom{2\ell}{\ell}$, diameter $D=\ell-1$, with the integer eigenvalues $\theta_i=(-1)^i(\ell-i)$ with corresponding multiplicities $m(\theta_i)=\binom{2\ell-1}{i}-\binom{2\ell-1}{i-1}$ for $i=0,1,\dots, \ell-1$. We note that $O_2$ is a triangle while $O_3$ is the Petersen graph. The girth of $O_\ell$ is $3$ if $\ell = 2$, 
$5$ if $\ell = 3$ and 
$6$ if $\ell > 3$. See, for instance, Biggs \cite{Biggs} and Godsil \cite{Godsil} for more details. Observe that, for $\ell$ odd, the distinct eigenvalues of $O_{\ell}$ are $$\ev(O_{\ell})=\{\ell, -(\ell-1),(\ell-2), \dots, 3,-2,1 \}=\{\ell, \ell-2,\dots, 3,1,-2,\dots, 1-\ell \},$$
while for $\ell$ even, the distinct eigenvalues of $O_{\ell}$ are $$\ev(O_{\ell})=\{\ell, -(\ell-1),(\ell-2), \dots, -3,2,-1 \}=\{\ell, \ell-2,\dots, 2,-1,-3,\dots, 1-\ell \}.$$
For $\ell>2$, we have that $\triangle=0$. Hence, for odd $\ell>2$, we have $\theta_s=1$, and $\theta_{s+1}=-2$. Similarly, for even   $\ell>2$, we have $\theta_s=2$, and $\theta_{s+1}=-1$.
Thus, by Theorem \ref{Gmaink3}, we have that, for $\ell$ odd,
\begin{align}
\alpha_3(O_\ell) \le  \binom{2\ell}{\ell} \frac{\ell^2-2\ell+2}{2(\ell+2)(\ell-1)(2\ell-1)} \label{odd1}
\end{align}
and, for $\ell$ even,
\begin{align}
\alpha_3(O_\ell) \le  \binom{2\ell}{\ell} \frac{\ell^2-4\ell+2}{2(\ell+1)(\ell-2)(2\ell-1)}.\label{odd2}
\end{align} 

\subsection{Johnson graphs}
Recall that the Johnson graph $J(n,k)$ has as vertices the $k$-element subsets of an $n$-element set, where two vertices are adjacent whenever the corresponding subsets intersect in exactly $k - 1$ elements. It has diameter $d=\min \{k,n-k\}$ and eigenvalues $\theta_j=(k-j)(n-k-j)-j$ for $0\le j\le d$ (Brouwer, Cohen and Neumaier~\cite{brouwer1989}). So, we have $\theta_{0}=k(k-k)$, and $\theta_d=-k$ if $d=k$ and $\theta_d=k-n$ if $d=n-k$. Each vertex $v$ of $J(n,k)$ is contained in $\frac{1}{2}k(n-k)(n-2)$ triangles, hence, $\triangle=k(n-k)(n-2)$. Thus, $\theta_{s}$ is the least eigenvalue such that 
$$\theta_{s}\ge 
\frac{k(n-k)-(n+d)+2}{d-1}.
$$
We can find the exact value of $s$, it is the floor of the smallest root of the quadratic equation $(k-j)(n-k-j)-j= \frac{k(n-k)-(n+d)+2}{d-1}$ in $j$.
In particular, for $J(2k,k)$, that is, if $n=2k$,  then $\theta_{s}$ is the least eigenvalue such that $\theta_s\ge k-2$ and the index $s$ is $s=\lfloor\frac{2k+1-\sqrt{8k-7}}{2}\rfloor$. Hence, by Theorem \ref{Gmaink3}, the Johnson graph $J(2k,k)$ satisfies
\begin{align}
\alpha_3(J(2k,k))&\le \frac{1}{k+1}\binom{2k}{k}\frac{\theta_s\theta_{s+1}-(\theta_s+\theta_{s+1}+2-3k)}{(k^2-\theta_s)(k^2-\theta_{s+1})}, \label{J2k}
%\\&=\frac{1}{k+1}\binom{2k}{k}\frac{\theta_s\theta_{s+1}-k(2k^2+2s^2-4ks-5k+2)}{(2ks-s^2+s)(2ks+2k-s^2-s)}
\end{align}
where $s=\lfloor\frac{2k+1-\sqrt{8k-7}}{2}\rfloor$ and $\theta_s=(k-s)^2-s.$
%\rm{\begin{example} \label{ex} 

In Table~\ref{tableegJ}, we show the bounds of the $3$-independence numbers for the Johnson graph $J(2k,k)$ for $k=2, 3,\dots, 7$ using Equation \ref{J2k}. 
	\begin{table}[h!]
	\caption{Bound for some Johnson graphs} 
	\centering
	{\begin{tabular}{|c|c|c|c|c|}
			\hline 
			Johnson Graph $J(2k,k)$ & $\theta_s$ & $\theta_{s+1}$ &Theorem~\ref{Gmaink3} & $\alpha_3$\\ 
			\hline 
			J(4,2) & 0 & -2 &  1 & 1 \\ 
			\hline 
			J(6,3) & 3 & -1 & 1 &1 \\ 
			\hline 
			J(8,4) & 2 & -1 & 2 &2\\ 
			\hline 
			J(10,5) & 7 & 1 & 3.11 & 2 \\ 
			\hline 
			J(12,6) & 6& 0 & 7.33 & 4 \\ 
			\hline 
			J(14,7) & 5 & -1 & 19.5& 8\\ 
		%	\hline 
		%	J(16,8) & 12 & 4 & 44 &-\\ 
		%	\hline 
		%	J(18,9) & 11 & 3 & 117.54 &-\\
		%	\hline
		%	J(20,10) & 10 & 2 & 342.78 &-\\
			\hline
	\end{tabular} }\label{tableegJ}
\end{table} 
%\end{example}}

%\begin{example}
%	We illustrate Theorem~\ref{Gmaink3} for some graphs in Table~\ref{tableeg}. We determine the value of $\theta_s$, and give the values of $b$ and $c$ of the polynomial, and then the bound $B(p)$.
%	%!htbp
%	\begin{table}[h!]
%		\caption{Bounds for some graphs} 
%		\centering
%		{\begin{tabular}{|c|c|c|c|c|}
%				
%				\hline 
%				Graph & $\theta_s$ & b & c & Bound \\ 
%				\hline 
%				J(7,2) & 3 & 1 & 7 & 1 \\ 
%				\hline 
%				J(8,4) & 2 & 4 & -4 & 2 \\ 
%				\hline 
%				J(14,7) & 5 & 3 & -33 & 19\\ 
%				\hline 
%				Cube(8) & 0 & 10 & 16 & 16 \\ 
%				\hline 
%				Cell-120 & 1 & 1.3 & -3.9 & 57 \\ 
%				\hline 
%				Ljubljana & 0 & 3.9 & 2.6 & 18\\ 
%				\hline 
%				OddGraph(5) & 1 & 5 & 2 & 8 \\ 
%				\hline 
%		\end{tabular} }\label{tableeg}
%	\end{table} 
%\end{example}
In Table~\ref{tablecom}, we compare our bounds with the bounds in Abiad, Coutinho, Fiol \cite{abiad2019k} and  Fiol~\cite{fiol2}, and the actual  $3$-independence number (when known) for some named graphs. %In the table, a period (.) after a number indicates that the bound is not an integer and has been rounded down.
 
\fontsize{9.5}{12}\selectfont
\begin{longtable}{|c|c|c|c|c|}
	\caption{Comparison between different bounds for the 3-independence number of some graphs.}\label{tablecom}\\
	\hline 
	Graph	& Bound in (\ref{genk}) \cite{abiad2019k} & Bound in (\ref{fiolk33}) \cite{fiol2}    & Theorem~\ref{Gmaink3} & $\alpha_3$ \\ 
	\hline 
	Johnson Graph J(14,7)	& 26.74  & 19.5  & 19.5 & 8\\ 
	\hline 
	Cube Graph(8)	& 114.25 & 16 & 16 &  16\\ 
	\hline 
	Odd Graph(6)	& 141.27 & 21  & 21 & 15 \\ 
	\hline 
	Balaban 10-cage	& 28 & 12.82 & 11.67 &  9\\ 
	\hline  
	Frucht graph & 3.62 & 2.35 & 2.25 & 2 \\ 
	\hline 
	Meredith graph & 27.48 & 9.18  &  8.58 &  7\\ 
	\hline 
	Moebius-Kantor graph & 6.4 & 2 & 2 &  2\\ 
	\hline  
	Bidiakis cube & 3.66 &1.92  & 1.50 & 1 \\ 
	\hline 
	Gosset graph & 1.31 & 1  & 1 &  1 \\ 
	\hline  
	Balaban 11-cage & 41.25 & 20.44 & 18.04 &  16\\ 
	\hline  
	Gray graph & 21.6 & 9 & 9 & 9 \\ 
	\hline 
	Nauru graph & 9.6  & 4  & 4 &  4\\ 
	\hline 
	Blanusa first snark graph & 5.68 & 2.43 & 2.43 &  2\\ 
	\hline 
	Pappus graph &7.2  & 3  & 3 &  3\\ 
	\hline 
	Blanusa second snark graph &4.82 & 3.15  & 2.50 & 2 \\ 
	\hline 
	Brinkmann graph & 4.49 & 2.14  & 1.97 & 1  \\ 
	\hline 
	Harborth graph & 12.83 & 8.49 & 8.11 &  6\\ 
	\hline 
	Perkel graph & 5.52 & 3.79  & 1 &  1\\ 
	\hline 
	Harries graph & 28 & 12.82  & 10.71 &  10\\ 
	\hline 
	Bucky ball & 19.02 &10.42 & 8.84 & 7 \\ 
	\hline 
	Harries-Wong graph & 28 &12.82 & 10.71 &  9\\ 
	\hline 
	Robertson graph & 3.91 & 1.82  & 1.77 & 1 \\ 
	\hline 
	Heawood graph & 5.37 &1  & 1 & 1 \\ 
	\hline 
	Cell 600 & 7.04 & 5.84  & 4.85 & 3 \\ 
	\hline 
	Cell 120 & 129.22 &74.51 & 57.41 & - \\ 
	\hline 
	Hoffman graph & 6.59 & 2 & 2 &  2\\ 
	\hline 
	Sylvester graph & 5.32 & 1 & 1 & 1 \\ 
	\hline 
	Coxeter graph & 7.70 & 4 & 4 & 4 \\ 
	\hline 
	Holt graph & 4.83 & 2.85 & 2.21 & 1 \\ 
	\hline 
	Szekeres snark graph & 15.76 & 8.46 & 7.35 &  6\\ 
	\hline 
	Desargues graph & 8 & 2.5 & 2.5 & 2 \\ 
	\hline 
	Horton graph & 38.4 & 16.89 & 16 &  14\\ 
	\hline 
	Dejter graph & 48.43 & 8 & 8 & 8 \\ 
	\hline 
	Tietze graph & 3.54 & 1.87 & 1.87 & 1 \\ 
	\hline 
	Double star snark & 9.29 &5.36 & 4.37 & 4 \\ 
	\hline 
	Durer graph & 3.40 & 2.24  & 2.18 &  2\\ 
	\hline 
	Klein 3-regular Graph & 17.09 & 9.80  & 7.29 &  7\\ 
	\hline 
	Truncated tetrahedron & 2.93 &2.2 & 1.6 &  1\\ 
	\hline 
	Dyck graph & 12.8 & 4 & 4 &  4\\ 
	\hline 
	Klein 7-regular graph & 2.04 & 1 & 1 & 1 \\ 
	\hline 
	Tutte 12-cage & 50.4 & 21 & 21 & 21  \\ 
	\hline 
	Ellingham-Horton 54-graph & 21.6  & 9.92  & 9 & 8 \\ 
	\hline 
	Tutte-Coxeter graph &12  & 5 & 5 &  5\\ 
	\hline 
	Ellingham-Horton 78-graph &31.2 & 13.73  & 13 &  11\\ 
	\hline 
	Ljubljana graph & 44.8 &21.37 & 18.67 &  17\\ 
	\hline 
	Tutte graph & 15.21 &8.17 & 7.17 & 6 \\ 
	\hline 
	F26A graph & 10.4 &5.10 & 3.61 & 3 \\ 
	\hline 
	Watkins snark graph & 14.75 &9.23  & 7.25 &  6\\ 
	\hline 
	Flower snark & 6.78 &3.55 & 2.89 &  2\\ 
	\hline 
	Markstroem graph & 6.64 &4.76 & 4.30 & 3 \\ 
	\hline 
	Wells graph & 4.72  & 2  & 2 &  2\\ 
	\hline 
	Folkman graph & 8.24 & 2.5 & 2.5 &  2\\ 
	\hline 
	Foster graph & 36 & 15 & 15 &  15\\ 
	\hline 
	McGee graph & 7.32 &4.13 & 3.0 &  2\\ 
	\hline 
	Franklin graph & 4.8 & 1.5 & 1.5 & 1 \\ 
	\hline 
	Hexahedron & 3.2 & 1  & 1 & 1 \\ 
	\hline 
	Dodecahedron & 4.52 &3.24 & 2.36 & 2 \\ 
	\hline 
	Icosahedron& 1.72 &1 & 1 & 1 \\ 
	\hline 
	Brouwer-Haemers Graph & 3.17 &1 & 1 & 1 \\ 
	\hline 
\end{longtable}
\normalsize

From Table~\ref{tablecom}, we observe that our bound holds with equality for some graphs, which
gives rise to a regular partition of the matrix  $A^3-(\theta_s+\theta_{s+1}+\theta_d)A^2+(\theta_d\theta_s+\theta_d\theta_{s+1}+\theta_s\theta_{s+1})A$, where $A$ is the adjacency matrix of the graphs. The graphs are listed in Table~\ref{regparts}.
 
\begin{table}[h!]
\centering
\caption{Some named graphs in Sage whose adjacency matrix has a regular partition.}\label{regparts}
%\fontsize{12}{12}\selectfont
\begin{tabular}{|c|c|}
	
	\hline 
	Graph	& $\alpha_3$ \\ 
	\hline 
	Cube Graph(8)	&  16\\ 
	\hline 
	Moebius-Kantor graph &  2\\ 
	\hline  
	Gosset graph &  1 \\ 
	\hline  
	Gray graph &  9 \\ 
	\hline 
	Nauru graph &   4\\ 
	\hline 
	Pappus graph & 3\\ 
	\hline 
	Perkel graph & 1\\ 
	\hline 
	Heawood graph & 1 \\ 
	\hline 
	Hoffman graph &   2\\ 
	\hline 
	Sylvester graph & 1 \\ 
	\hline 
	Coxeter graph &  4 \\ 
	\hline 
	Dejter graph &  8 \\ 
	\hline 
	Dyck graph &   4\\ 
	\hline 
	Tutte 12-cage &  21  \\ 
	\hline 
	Tutte-Coxeter graph & 5\\ 
	\hline 
	Wells graph &  2\\ 
	\hline 
	Foster graph &  15\\ 
	\hline 
	Hexahedron & 1 \\ 
	\hline 
	Icosahedron& 1 \\ 
	\hline 
	Brouwer-Haemers Graph & 1 \\ 
	\hline 
\end{tabular}
\end{table}
\normalsize

%\begin{eg} \label{antipodal1}
%	Let $G$ be an antipodal bipartite distance-regular graph, with degree  $d$ and diameter $3$. These graphs have $n=2(d+1)$ vertices \cite{brouwer1989}, intersection array $\{d,d-1,1;1,d-1,d\}$, and distinct eigenvalues \begin{align}
%	\theta_{0}=d, \, \theta_{1}=1, \, \theta_{2}=-1, \, \theta_{3}=-d.
%	\end{align}
\section{Comparing our bound with the bound of Fiol}\label{fiolcompare}
Recall the polynomial given by Fiol~\cite{fiol2} and its corresponding bound on the $3$-independence number presented in Corollary~\ref{fiolk33}. Note that $\theta_i$ in Corollary~\ref{fiolk33} is the largest eigenvalue such that $\theta_i\leq -1$. Recall also that from Theorem~\ref{Gmaink3}, $s$ is the largest index such that $\theta_s\ge -\frac{\theta_{0}^2+\theta_{0}\theta_d-\triangle}{\theta_{0}(\theta_d+1)}$, where $\triangle=\text{max}_{u\in V}\{(A^3)_{uu}\}$. %, that is, twice the largest number of triangles on a vertex in $G$.  
Let $\tau= -\frac{\theta_{0}^2+\theta_{0}\theta_d-\triangle}{\theta_{0}(\theta_d+1)}$. Note that $\triangle$ is $2n_t$ for a graph $G$ that is at least $3$-partially walk regular.  Since $\triangle\le \theta_{0}(\theta_{0}-1)$, we have that $\theta_s\ge \tau \ge -1$. In particular, if $G$ is bipartite, then $\theta_s\ge 0$.  If  there exists an eigenvalue $\mu$ of $G$, where $\mu\in (\theta_i,\theta_s)$, then $\theta_{s+1}\ge \mu$ and $\theta_{i-1}\le \mu$. Thus, $\theta_{i-1}\ne \theta_s$ (that is, $\theta_i\ne \theta_{s+1}$), and so these two bounds, Corollary~\ref{fiolk33} and Theorem~\ref{Gmaink3}, do not coincide. 

Let us summarize the above discussion in the following proposition.
\begin{proposition} Let $\tau= -\frac{\theta_{0}^2+\theta_{0}\theta_d-\triangle}{\theta_{0}(\theta_d+1)}$. The bound on the $3$-independence number given in Corollary~\ref{fiolk33} is greater than the bound in Theorem~\ref{Gmaink3} if and only if there exists an eigenvalue $\mu$ of $G$ such that $\mu\in (-1,\tau)$.
\end{proposition}
Alternatively, if we take the polynomial $f(x)$ in Fiol's bound and take $\theta_i$ to be $\theta_{s+1}$, then the polynomial is of the form $f(x)=B\,p(x)+C$ for some constants $B$ and $C$, and $p(x)$ is  our polynomial. %And, $B=1/((\theta_0-\theta_d)(\theta_0-\theta_i)(\theta_0-\theta_{i-1}))$ and $C=-\theta_i\theta_{i-1}\theta_d/((\theta_0-\theta_d)(\theta_0-\theta_i)(\theta_0-\theta_{i-1}))$. 
So, the bounds resulting from $f(x)$ and $p(x)$ will coincide precisely when $\theta_i$ coincides with $\theta_{s+1}$. %Thus our index is the most accurate choice. %That is, the best bound in Theorem \ref{fiolsgen} is obtained by choosing the index in $I=\{s,s+1,d\}$.

\section{Acknowledgments}
The authors acknowledge helpful comments
and suggestions of  Mateja \v{S}ajna. 
We appreciate the anonymous referees for their valuable comments and suggestions, improving the presentation of this paper.
Research is supported by NSERC Canada. 

\bibliographystyle{plain}

\bibliography{references}

\end{document}